\documentclass[10pt,journal,compsoc]{IEEEtran}

\usepackage{amsmath,bm, amsfonts,color,colordvi,caption,subcaption,url}
\usepackage{graphicx}
\usepackage{mathrsfs}
\usepackage{CJK}
 \usepackage{amssymb}
\usepackage{makecell}
\usepackage{picinpar}
\usepackage{multirow}
\usepackage{picinpar}

\usepackage{algorithm}
\usepackage{algorithmic}
\usepackage{ragged2e}
\usepackage{array}
\RequirePackage{xcolor}[2.11]
\colorlet{siaminlinkcolor}{green!50!black}
\colorlet{siamexlinkcolor}{red!60!black}
\colorlet{siamreviewcolor}{black!50}

\usepackage[colorlinks=true,
    linkcolor=cyan,
    citecolor=siamexlinkcolor,      
    urlcolor=cyan]
    {hyperref}
\usepackage[noadjust]{cite}
\newcounter{RomanNumber}

\makeatletter

\newcommand{\Rmnum}[1]{\expandafter\@slowromancap\romannumeral #1@}
\makeatother
\graphicspath{ {./images/} }

\newtheorem{theorem}{Theorem}[section]

\newtheorem{lemma}{Lemma}[section]

\newtheorem{definition}{Definition}[section]
\newtheorem{remark}{Remark}[section]
\newtheorem{example}{Example}[section]
\def\T{{\mathbb T}}
\def\R{{\mathbb R}}
\def\P{{\mathbb P}}
\def\bfe{{\bf 1}}

\def\bfd{{\bf d}}
\def\bfx{{\bf x}}
\def\bfw{{\bf w}}
\def\bfy{{\bf y}}
\def\bfz{{\bf z}}
\def\bfu{{\bf u}}
\def\bfal{{\bf {\boldsymbol \alpha}}}
\def\supp{{\rm supp}}

\usepackage[nameinlink,noabbrev,capitalize]{cleveref}

\ifCLASSINFOpdf
\else
\fi

\ifCLASSOPTIONcompsoc
\else
  \usepackage{cite}
\fi

\begin{document}
\title{Sparse SVM for Sufficient Data Reduction}

\author{
{Shenglong Zhou}% <-this % stops a space
\IEEEcompsocitemizethanks{
\IEEEcompsocthanksitem S.L. Zhou is with the Department of Electrical and Electronic Engineering, Imperial College London, London, UK. Email: slzhou2021@163.com.%\protect
%\IEEEcompsocthanksitem * Corresponding author
}

\thanks{Manuscript received 25 Aug 2020; revised 4 Feb 2021; accepted 21 April 2021.}}

\IEEEtitleabstractindextext{\justify
\begin{abstract}
Kernel-based methods for support vector machines (SVM) have shown highly advantageous performance in various applications. However, they may incur prohibitive computational costs for large-scale sample datasets. Therefore, data reduction (reducing the number of support vectors) appears to be necessary, which gives rise to the topic of the sparse SVM. Motivated by this problem, the sparsity constrained kernel SVM optimization has been considered in this paper in order to control the number of support vectors. Based on the established optimality conditions associated with the stationary equations, a Newton-type method is developed to handle the sparsity constrained optimization. This method is found to enjoy  the one-step convergence property if the starting point is chosen to be close to a local region of a stationary point, thereby leading to a super-high computational speed. Numerical comparisons with several powerful solvers demonstrate that the proposed method performs exceptionally well, particularly for large-scale datasets in terms of a much lower number of support vectors and shorter computational time.
\end{abstract}

% Note that keywords are not normally used for peer review papers.
}

\maketitle


% is not selected <OR> if conference mode is selected - because compsoc
% conference papers position the abstract like regular (non-compsoc)
\IEEEdisplaynontitleabstractindextext

% \fi
\IEEEpeerreviewmaketitle

\ifCLASSOPTIONcompsoc
\IEEEraisesectionheading{\section{Introduction}\label{sec:introduction}}
\else
\section{Introduction}
\label{sec:introduction}
\fi

\IEEEPARstart{S}{upport} vector machines (SVM) were first introduced by Cortes  and Vapnik \cite{CV95} and are currently popular classification tools in machine learning, statistic and pattern recognition. The goal is to find a hyperplane in the input space that best separates the training dataset to enable accurate prediction of the class of some newly input data. The paper focuses on the binary classification problem:  suppose that we are given a training dataset $\{(\bfx_i,y_i):i=1,2,\cdots,m\},$ where $\bfx_i\in \R^{n}$ is the sample vector and $y_i\in \{-1,1\}$ is the binary class.  The task is to train a hyperplane $\langle \bfw, \bfx\rangle+b=w_1x_1+\cdots+w_nx_n+b=0$ with variable $\bfw\in \R^{n}$ and bias $b\in\R$ to be estimated based on the training dataset. For any newly input vector $\overline\bfx $, one can predict the corresponding class $\overline y $ by $\overline y =1$ if $\langle \bfw, \overline \bfx \rangle+b>0$ and $\overline y =-1$ otherwise.  There are two possible scenarios to find an optimal hyperplane: linearly separable and  inseparable training datasets in the input space.  For the latter, the popular approach is to solve the so-called soft-margin SVM optimization.

There is a vast body of work on the design of the loss functions $\ell$ for dealing with the soft-margin SVM optimization. One of the most well-known loss functions is the Hinge loss, giving rise to the Hinge soft-margin SVM model. To address such a problem, the dual kernel-based SVM optimization is usually used, for which the objective function involves an $m\times m$-order Gram matrix known as the kernel matrix. The complexity of computing this kernel matrix is approximately $O(m^2n)$, thereby making it impossible to develop methods for training on  million-size data since data  storage on such a scale requires  large memory and the incurred computational cost is prohibitively high. Therefore, to overcome this drawback, data reduction has drawn much attention with the goal of using  a small portion of samples to train a classifier.

It is well known that by the Representer Theorem, the  classifier $\bfw^*$ can be expressed as \begin{eqnarray}
\label{Representer-Theorem}\bfw^*=\sum_{i=1}^m \alpha_i^* y_i\bfx_{i},  \end{eqnarray}
 where $\bfal^*$ is a solution to the dual kernel-based  SVM optimization. The training vectors $\bfx_{i}$
corresponding to nonzero $\alpha_i^*$ are known as the support vectors. If large numbers of coefficients $\alpha_i^*$ are zeros, then the number of the support vectors can be reduced significantly. Therefore, computations and storage for large scale size data are possible since only the support vectors are used. However, solutions to the dual kernel-based SVM optimization are not sparse enough generally.  An impressive number of approaches have been developed for guaranteeing that the solution is sufficiently sparse. The methods aiming to reduce the number of support vectors can be categorized as the sparse SVM group.  We will explore more in the sequel. 

\subsection{Selective Literature Review}
Extensive work has focused on designing loss functions  $\ell$ to cast   efficient soft-margin SVM models that can be summarized into two classes based on the convexity of $\ell$. Convex soft margin loss functions include the famous hinge loss \cite{CV95}, the pinball loss \cite{JHS2013,HSS2014}, the hybrid Huber loss \cite{RZ2007,WZZ2008, XAC2016},  the square loss \cite{SV1999, YTH2014}, the exponential loss \cite{FS1997} and log loss \cite{FHT2000}. Convexity  makes the computations of their corresponding SVM models tractable, but it also induces the unboundedness, thereby reducing the robustness of these functions to the outliers from the training data.  To overcome this drawback, \cite{MBB2000},\cite{PNF2003} set an upper bound and enforce the loss functions to stop increasing after a certain point. This makes the convex loss functions become nonconvex. Other nonconvex losses include the ramp loss  \cite{CFWB2006}, the {truncated pinball loss} \cite{SNZT2017}, the {asymmetrical truncated pinball loss}  \cite{YD2018}, the {sigmoid loss}  \cite{PNP2000}, and the normalized sigmoid cost loss  \cite{MBBF1999}.  Compared with the convex margin loss functions, most  nonconvex functions are less sensitive to the outliers due to their boundedness. However, nonconvexity generally gives rise to difficulties in numerical computations.

A research effort following a different direction investigated methods for data reduction, namely, reducing the number of the support vectors, giving rise to the topic of the sparse SVM. One of the earliest attempts can be traced back to \cite{OF98}, where it was suggested to solve the kernel SVM optimization problem to find a solution first and  then seek a sparse approximation through support vector regression. This idea was then adopted as a key component of the method developed in \cite{ZS05}.% where the training method can exclude the samples that incurred the separation hypersurface highly convoluted so that a small number of the support vectors were sufficient for  describing a less convoluted hypersurface for separating two classes. 
The famous reduced SVM (RSVM)  in \cite{LM01}  randomly picked a subset of the training set,
and then searched for a solution (supported only on the
picked training set) to a smooth SVM optimization that minimized the loss on the entire training set. In \cite{BBEB03}, from \eqref{Representer-Theorem},  they substituted $\bfw$ by the expression $\sum_{i=1}^m \alpha_i y_i\bfx_{i}$   and employed
the $\ell_1$-norm regularization  of  $\bfal$ in the soft-margin loss model, leading to effective variable/sample selections because the $\ell_1$-norm regularization can render a sparse structure of the solution \cite{T96}. Similarly,  \cite{KCD06} also took advantage of the expression of  \eqref{Representer-Theorem}. They performed a greedy method, where in each step,  a new training sample was carefully selected into the set of the training vectors to form a new subproblem. Since the number of the training vectors was relatively small in comparison with the total size of the training samples,  the subproblem was on small scale  and thus can be addressed by Newton method quickly.  In \cite{CSS13}, a subgradient descent algorithm was proposed where in each step, only the samples with the correct classification inside the margin but maximizing a gap were selected.  Other relevant methods include the so-called   reduced set  methods  \cite{B96,BS97}, the Forgetron algorithm \cite{DSS06},  the condensed vector machines training method \cite{NMTH10} and those in \cite{KH04,WSB05}. Numerical experiments have demonstrated that these methods perform exceptionally well for the reduction in the number of support vectors.

 {We note that some sparse SVM methods \cite{LM01,KCD06} aim to reduce the size of the dual problem, a quadratic kernel-based optimization,  to reduce computational cost and the required storage. For the same purpose, two alternatives are available for dealing with data on large scales. The first is to process the data prior to employing a method. For instance,  \cite{WLLZZ14, PCW06}   select informative samples and remove the useless samples to reduce the computational cost but while preserving the accuracy. For another example, the `cascade' SVM \cite{GCBDV04} and  the mixture SVM  \cite{CBB02} split the dataset into several disjointed subsets and then carry out parallel optimization to accelerate the computation.  The second approach effectively benefits from kernel tricks. Nystr$\ddot{\rm o}$m method  \cite{WS01} and the sparse greedy approximation technique \cite{SS00} can be adopted to construct the Gram matrix to reduce the cost of solving the kernel-based dual problem. Moreover, \cite{FS01,ZSW20} exploited the low-rank kernel representations to make the interior point method tractable for large scale datasets.
}
\subsection{Methodology}
Mathematically, the soft-margin SVM takes the form of
\begin{eqnarray}\label{SM-SVM}
\min_{\bfw\in \R^{n}, {b\in\R}}~~  \frac{1}{2}\|\bfw\|^2+C\sum_{i=1}^{m}\ell\Big[1-y_{i} (\langle \bfw, \bfx_i\rangle +b) \Big],
\end{eqnarray}
 where $C>0$ is a penalty parameter, $\|\cdot\|$ is the Euclidean norm and $\ell$ is a loss function. One of the most famous loss functions is the Hinge loss $\ell_{\rm H}(t):=\max\{0, t\}$, giving rise to the Hinge loss soft-margin SVM, for which the dual problem is the following quadratic kernel-based SVM optimization,
\begin{eqnarray}
\label{SM-SVM-hinge-dual}
\min_{\bfal\in \R^{m}}&& d(\bfal):=\frac{1}{2} \sum_{i=1}^m\sum_{j=1}^m\alpha_i\alpha_jy_iy_j\langle\bfx_i,\bfx_j\rangle-\sum_{i=1}^m\alpha_i,\nonumber\\
{\rm s.t.}&&\sum_{i=1}^m\alpha_iy_i=0,~~ 0\leq \alpha_i\leq C, i\in[m], 
\end{eqnarray} 
where $[m]:=\{1,\cdots,m\}$. In this paper, we focus on the following soft-margin SVM,
\begin{eqnarray}\label{SM-SVM-cc}
\min_{\bfw\in \R^{n}, {b\in\R}}~~  \frac{1}{2}\|\bfw\|^2+\sum_{i=1}^{m}\ell_{cC}\Big[1-y_{i} (\langle \bfw, \bfx_i\rangle +b) \Big],
\end{eqnarray}
where the soft-margin loss $\ell_{cC}$ is defined by
\begin{eqnarray}\ell_{cC}(t):= 
\begin{cases}
 {Ct^2}/{2},&t\geq0,\\
 ~{ct^2}/{2},&t<0.
\end{cases}
\end{eqnarray} 
Here, $c>0$ is chosen to be smaller than $C$, namely $c<C$. The soft-margin SVM  \eqref{SM-SVM-cc} implies that it gives penalty $C$ for $t_i:=1-y_{i} (\langle \bfw, \bfx_i\rangle +b)\geq0$ and $c$ otherwise.  Note that $\ell_{cC}$ is reduced to the squared Hinge loss $(\ell_{\rm H}(t))^2$  when $c=0$. \cref{the:duality} shows that  the dual problem of \eqref{SM-SVM-cc} takes the form of
\begin{eqnarray}
\label{SM-SVM-h-equ-dual}
\min_{\bfal\in \R^{m}} ~D(\bfal):=d(\bfal)+ \sum_{i=1}^{m}h_{cC}(\alpha_i),~~
{\rm s.t.}~\sum_{i=1}^m\alpha_iy_i=0,
\end{eqnarray}
where $h_{cC}$ is defined by
\begin{eqnarray}\label{hcc}
h_{cC}(t):= 
\begin{cases}
 {t^2}/{(2C)},&t\geq0,\\
 {t^2}/{(2c)},&t<0.
\end{cases}
\end{eqnarray} 
Motivated by the work in the sparse SVM \cite{OF98}-\cite{WSB05}, in this paper, we aim to solve the following sparsity constrained kernel-based SVM optimization problem,
\begin{eqnarray}
\label{SM-SVM-h-equ-dual-sparse}
\min_{\bfal\in \R^{m}}~ D(\bfal),~~{\rm s.t.}~\sum_{i=1}^m\alpha_iy_i=0, ~\|\bfal\|_0\leq s, 
\end{eqnarray}
where  $s\in[m]$ is a given integer satisfying $s\ll m$ and is called the sparsity level, and $\|\bfal\|_0$  is  the zero norm of $\bfal$, counting the number of nonzero elements of $\bfal$.   Compared with the classic kernel SVM model \eqref{SM-SVM-hinge-dual}, the problem \eqref{SM-SVM-h-equ-dual-sparse} has at least three advantages: 
A1) The objective function is strongly convex; hence, the optimal solution exists (see  \cref{existence-global-minimizers}) and  is unique  under mild conditions (see  \cref{nec-suff-opt-con-tau}). By contrast, the objective function of \eqref{SM-SVM-hinge-dual} is convex but not strongly convex when $n\leq m$. This means that it may have multiple optimal solutions.
 A2)    Since the bounded constraints $0\leq \alpha_i\leq C, i\in[m]$ are absent,  the computation  is somewhat more tractable. Note that when $m$ is large,  {these bounded constraints in  \eqref{SM-SVM-hinge-dual} may incur some computational costs if no additional efforts, such as the use of the perturbed KKT optimality conditions \cite{FS01,ZSW20}, are paid.}
 A3)   Most importantly, the constraint $\|\bfal\|_0\leq s$ means that at most $s$ nonzero elements are contained in $\bfal$, i.e., the number of the support vectors is expected to be less than $s$ by \eqref{Representer-Theorem}. Thus, the number can be controlled to be small.

\subsection{Contributions}
The contributions of this paper are summarized as follows.

C1) As mentioned above, the sparsity constrained model \eqref{SM-SVM-h-equ-dual-sparse} has its  advantages. To the best of our knowledge, this is the first paper that employs the sparsity constraint into the kernel SVM model. Such a constraint allows us to control the number of support vectors and hence to perform  sufficient data reduction, making extremely large scale computation possible and significantly reducing the demand for huge  hardware memory.

C2) Based on the established optimality condition associated with the stationary equations by  \cref{stationary-equation}, a Newton-type method is employed. The method shows very low computational complexity and enjoys one-step convergence property. This convergence result is much better than the locally quadratic convergence property that is typical of  Newton-type methods. Specifically, the method converges to a stationary point within one step if the chosen starting point is close enough to the stationary point, see  \cref{the:i-iterate-convergence}. Moreover,  since the sparsity level $s$ that is used to control  the number of support vectors is unknown in practice, the selection of $s\in[m]$ is somewhat tedious. However, we successfully design  a mechanism to tune the sparsity level $s$ adaptively.

C3) Numerical experiments have demonstrated that the proposed method performs exceptionally well particularly for datasets with $m\gg n$. Compared  with several leading solvers, it takes much shorter computational time due to a tiny number of support vectors used for large scale data.

\subsection{Preliminaries}
We present some notation to be employed throughout the paper in \cref{notation}.

\begin{table}[!th]
\renewcommand{\arraystretch}{1}\addtolength{\tabcolsep}{-4pt}
\begin{center}
\caption{List of notations}\label{notation}\vspace{-3mm}
\begin{tabular}{ll}\\\hline
Notation & Description\\\hline
$[m]$       & The index set $\{ 1,2,\cdots, m\}$.\\
$|T|$         & The number of elements of an index set $T\subseteq[m]$. \\
$ \overline T$        & The complementary set of an index set $T$, namely, $[m]\setminus T$. \\
$\bfal_{ T}$   & The sub-vector of $\bfal$ indexed on $T$ and $\bfal_{ T}\in\R^{|T|}.$\\
$|\bfal|$     & $:=(|\alpha|,\cdots,|\alpha|)^\top$. \\
%$\bfal_+$     & $:=((\alpha_1)_+,\cdots,(\alpha_m)_+)^\top$.\\
$\|\bfal\|_{[s]}$  & The $s$th largest  element of $|\bfal|$.\\
%$ \| \bfal \|_{[1]}$&  The $\ell_\infty$ norm of $\bfal$, namely, $ \max_i|\alpha_i|$.\\
%$\|\bfal\|_0$ & The number of non-zero elements of $\bfal$. \\
$\supp(\bfal)$& The support set of $\bfal$, namely, $\{i\in [m]: \alpha_i\neq0\}$.\\
$\bfy$       & The labels/classes $(y_{1},\ldots,y_{m})^{\top}\in \R^{m}$. \\
$X$       & The samples data $[\bfx_{1}~\cdots~\bfx_{m}]^\top\in \R^{m\times n}$. \\
$Q$          & $:=[y_1\bfx_{1}~\cdots~y_m\bfx_{m}]\in \R^{n\times m}$.\\
$Q_T$          &The sub-matrix containing the columns of $Q$ indexed on ${T}$.\\
$Q_{\Gamma,T}$  &The sub-matrix containing the rows of $Q_{T}$ indexed on ${\Gamma}$.\\
$I$       & The identity matrix.\\
$P$          &$:=Q^\top Q+I/C$.\\
$\bfe$       & $:=(1,\cdots,1)^{\top}$ whose dimension varies in the context.\\
%$ \| Q \|$& The maximum singular value of $Q$.\\
$N({\bfal},\delta)$& The neighbourhood of ${\bfal}$ with radius $\delta>0$, namely,\\
& $\{\bfu\in\R^m: \|\bfu-{\bfal}\|<\delta\}.$\\\hline
\end{tabular}\end{center}
\end{table}

\noindent Notation in \cref{notation} enables us to rewrite $D(\bfal)$ in \eqref{SM-SVM-h-equ-dual-sparse} as
\begin{eqnarray}
\label{obj}
D(\bfal)=\frac{1}{2}\|Q\bfal\|^2  +\frac{1}{2} \langle E(\bfal)\bfal,  \bfal \rangle - \langle\bf1,  \bfal\rangle,
\end{eqnarray}
where $E(\bfal)$ is a diagonal matrix with 
\begin{eqnarray}\label{diag-hessian}
E_{ii}(\bfal):=(E(\bfal))_{ii}= \left\{\begin{array}{lll}
 1/C,& \alpha_i\geq 0, \\
 1/c,& \alpha_i< 0.  
\end{array}\right.
\end{eqnarray}
The gradient $\nabla  D(\bfal)$ and Hessian $H(\bfal)$ of $D(\bfal)$ are
\begin{eqnarray}\label{gradient-hessian}
\arraycolsep=1.4pt\def\arraystretch{1.5}
\begin{array}{lll}
\nabla  D(\bfal)  :=  H(\bfal)\bfal-\bfe,&&
H(\bfal) := Q^\top Q + E(\bfal).\end{array}
\end{eqnarray}
It is observed that the Hessian matrix is positive definite for any $\bfal\in\R^m$ 
 and \begin{eqnarray}\label{hessian}
H  (\bfal)  \succeq Q^\top Q +I/C=P \succ 0,
\end{eqnarray}
due to $C>c$, where  $A\succeq0$ ($A\succ0$) means $A$ is semi-definite (definite) positive. Here, $A \succeq B$ represents that $A-B$ is semi-definite positive.  For notational convenience, hereafter, for a given $\bfal\in\R^m$ and $ b \in\R$, let
\begin{eqnarray}\label{z-z*-zk} 
\bfz&:=&(\bfal; b )=\left[ 
\begin{array}{c}
\bfal \\
 b 
\end{array}
 \right].
\end{eqnarray}
Based on this, we denote the following functions
\begin{eqnarray}\label{beta-mu} 
\arraycolsep=1.4pt\def\arraystretch{1.5}
\begin{array}{lll}
 g(\bfz)&:=& \nabla D(\bfal) +\bfy b \overset{\eqref{gradient-hessian}}{=} H(\bfal)\bfal-\bfe+\bfy b , \\
  g_T(\bfz)&:=&(g(\bfz))_T,~~~~
  H_T(\bfal):=(H(\bfal))_{TT}. \end{array}
\end{eqnarray}
Here, $g(\bfz)$ is a vector and $g_T(\bfz)$ is a sub-vector of $g(\bfz)$.  $H(\bfal)$ is a matrix and $H_T(\bfal)$ is the sub-principal matrix of $H(\bfal)$ indexed by $T$. Similar rules are also applied for $\bfz^*$ and $\bfz^k$.
\subsection{Organization}
The rest of the paper is organized as follows.   In the next section, we focus on the sparsity constrained model \eqref{SM-SVM-h-equ-dual-sparse}, establishing the optimality condition associated with the $\eta$-stationary point. The condition is then equivalently transferred to stationary equations that allow  us to adopt the Newton method in \cref{sec:SNM} where we also prove that the proposed method converges to  a stationary point within one step. Furthermore, a strategy of tuning the sparsity level $s$ is employed to derive {\tt NSSVM}. Numerical experiments are presented in  \cref{numerical}, where the
implementation of  {\tt NSSVM} as well as its comparisons with some leading solvers are provided. We draw the conclusion in the last section and present all of the proofs in the appendix.
\section{Optimality}\label{Sec:sparse}
Prior to the establishment of the optimality condition of the sparsity constrained kernel SVM optimization \eqref{SM-SVM-h-equ-dual-sparse}, 
hereafter, for a point $ \bfal^*\in\R^m$, its support set is always denoted by
\begin{eqnarray}
\label{SM-SVM-h-equ-dual-sparse-opt}
S_*&:=&\supp(\bfal^*).
\end{eqnarray} 
We now claim that \eqref{SM-SVM-h-equ-dual} is indeed the dual problem of \eqref{SM-SVM-cc}. 
\begin{theorem}\label{the:duality}
 Problem \eqref{SM-SVM-h-equ-dual} is the dual problem of \eqref{SM-SVM-cc} and admits a unique optimal solution, say $ \bfal^*$. Furthermore,  the optimal solution  of the primal problem  \eqref{SM-SVM-cc} is 
\begin{eqnarray} \label{w-b}
\widehat \bfw =Q \bfal^*,~~\widehat b = \frac{\langle \bfy_{S_*}, \bfe - H_{S_*}(\bfal^*)\bfal^*_{S_*}\rangle}{|S_*|}. 
\end{eqnarray}
\end{theorem}
Before we embark on the main theory in this section, we explore more of the sparse constraint. Let $\P_s$ be defined by
 \begin{eqnarray}\label{HTP-operator}
 \P_s(\bfal) ={\rm argmin}_{\bfu} \Big\{  \|\bfu-\bfal\|:   \|\bfu\|_0\leq s \Big\}, \end{eqnarray}
which can be obtained by retaining the $s$ largest elements in magnitude from $\bfal$ and  setting the remaining to zero. $\P_s$ is known as the Hard-thresholding operator. To well characterize the solution of \eqref{HTP-operator}, we define a useful set by
\begin{eqnarray}
\label{T-beta} 
\T_s(\bfal):=\left\{T\subseteq[m]: 
\begin{array}{ll}
|T|=s,\\ |\alpha_i|\geq|\alpha_j|,\forall i\in T,\forall j\notin T
\end{array}
 \right\}.
\end{eqnarray}
This definition of $\T_s$ allows us to express $\P_s$ as
\begin{eqnarray}
\label{P-s-T}~~~~~~~~
\P_s(\bfal):=\left\{\left(
\bfal_T; ~
0\right):~ T\in \T_s(\bfal)\right\}.
\end{eqnarray}
As the $s$th largest element  of $|\bfal|$ may not be unique,  $\T_s(\bfal)$ may have multiple elements, so does  $\P_s(\bfal)$. For example, consider  $\overline\bfal =(1,-1,0,0)^\top$.  It is easy to check that $\T_1(\overline\bfal)=\{\{1\},\{2\}\}$, $\P_1(\overline\bfal)=\{(1,0,0,0)^\top,(0,-1,0,0)^\top\}$ and
$\T_2(\overline\bfal)=\{\{1,2\}\}$, $\P_2(\overline\bfal)=\{\overline\bfal\}$. 
\subsection{$\eta$-Stationary Point}
In this part, we focus on the sparsity constrained kernel SVM problem \eqref{SM-SVM-h-equ-dual-sparse}. First, we  conclude that it admits a global solution/minimizer.
\begin{theorem}\label{existence-global-minimizers} The global minimizers of \eqref{SM-SVM-h-equ-dual-sparse} exist. 
\end{theorem}
{The Lagrangian function of \eqref{SM-SVM-h-equ-dual-sparse} is $L(\bfal,b):=D(\bfal) + b\langle\bfal,\bfy\rangle$, where $b$ is the Lagrangian multiplier. To find a solution to  \eqref{SM-SVM-h-equ-dual-sparse}, we consider its
Lagrangian dual problem 
\begin{eqnarray}
\label{SM-SVM-h-equ-dual-sparse-lag}
\max_{b\in\R} \min_{\bfal\in \R^{m}} \{ L(\bfal,b): \|\bfal\|_0\leq s\}. 
\end{eqnarray}
We note that $\nabla_{\bfal}L(\bfal,b)= g(\bfz)$ from \eqref{beta-mu} and $\nabla_{b}L(\bfal,b)=\langle\bfal,\bfy\rangle$. Based on these, 
we introduce the concept of an $\eta$-stationary point of the problem \eqref{SM-SVM-h-equ-dual-sparse}. }
 \begin{definition}
 We say $\bfal^*$ is an $\eta$-stationary point of \eqref{SM-SVM-h-equ-dual-sparse} with $\eta>0$ if there is  $ b^*\in\R$ such that
\begin{eqnarray}
\label{SM-SVM-h-equ-dual-sparse-tau}\left\{\begin{array}{cll}
\bfal^*&\in&\P_s\Big[\bfal^*-\eta g( \bfz^*) \Big],\\ 
0&=&\langle\bfal^*, \bfy \rangle. 
 \end{array}\right.
\end{eqnarray}
 \end{definition}
From the notation \eqref{z-z*-zk}  that $
\bfz^*=(\bfal^*; b ^*)$. We also say  $\bfz^*$ is an $\eta$-stationary point of \eqref{SM-SVM-h-equ-dual-sparse} if it satisfies the above conditions. We characterize an $\eta$-stationary point   of \eqref{SM-SVM-h-equ-dual-sparse} as a system of equations through the following theorem.
\begin{theorem}\label{stationary-equation}A point $\bfz^*$ is an $\eta$-stationary point  of \eqref{SM-SVM-h-equ-dual-sparse} with $\eta>0$ if and only if $\exists~T_*\in \T_s(\bfal^*-\eta g ( \bfz^*) )$ such that 
\begin{eqnarray}
\label{stationary-equation-T}
F(\bfz^*;T_*)&:=&\left[
\begin{array}{c}
  g_{T_*}(\bfz^*)\\ 
  \bfal_{\overline T_*}^*  \\
   \langle\bfal_{T_*}^*, \bfy_{T_*} \rangle
\end{array}\right]=0.
\end{eqnarray}
\end{theorem}
We call \eqref{stationary-equation-T} the stationary equations that allow for employing the Newton method  { used to find solutions to system of equations. It is well known that the (generalized) Jacobian of the involved equations is required to update the Newton direction. Therefore, comparing with the conditions in \eqref{SM-SVM-h-equ-dual-sparse-tau}, equations in \eqref{stationary-equation-T} enable easier implementation of the Newton method.  More specifically,   \eqref{SM-SVM-h-equ-dual-sparse-tau} involves inclusions and  its Jacobian is difficult to obtain as $\P_s(\cdot)$ is nondifferentiable. By contrast, for any fixed $T\in \T_s(\bfal-\eta g ( \bfz) )$, all functions in $F(\bfz;T)$ are  differentiable and the Jacobian matrix of $F(\bfz;T)$ can be derived by,
\begin{eqnarray}
\label{Jacobian-F}
\nabla F(\bfz;T)=\left[
\begin{array}{ccc}
 H_{T}(\bfal) & 0 & \bfy_T\\ 
0&I&0\\
  \bfy_T^\top&0&0  
\end{array}\right],
\end{eqnarray}
which is always nonsingular due to $H_{T}(\bfal)\succ 0$.}
\begin{remark} What is the role of $b^*$? The second condition in \eqref{stationary-equation-T}  indicates $S_*\subseteq T_*.$ Then the first equation in \eqref{stationary-equation-T} and \eqref{beta-mu} yields that $$
0=g_{S_*}( \bfz^*)=H_{S_*}(\bfal^*)\bfal^*_{S_*}-\bfe+\bfy_{S_*} b ^*,$$ 
which together with \eqref{w-b} means $b^*=\widehat{b}$. Namely,  $b^*$ is the bias. 
Therefore, seeking an $\eta$-stationary point can acquire the solution to the dual problem \eqref{SM-SVM-h-equ-dual-sparse} and the bias $b^*$ to the primal problem \eqref{SM-SVM-cc} at the same time. 
\end{remark}
We now reveal the relationships among local minimizers,  global minimizers and  $\eta$-stationary points. These relationships indicate that to pursue a local (or even a global) minimizer,  {we instead find an $\eta$-stationary point because it could be effectively derived  from \eqref{stationary-equation-T} by the numerical method proposed in the next section.}
\begin{theorem}\label{nec-suff-opt-con-tau} Consider a feasible point $\bfal^*$ to the problem \eqref{SM-SVM-h-equ-dual-sparse}, namely, $\|\bfal^*\|_0\leq s$ and $\langle\bfal^*, \bfy \rangle=0$. We have the following relationships. 

\noindent i) For local minimizers and $\eta$-stationary points, it has
\begin{itemize}
\item[$a)$] an $\eta$-stationary point $\bfal^*$ is also  a local minimizer.
\item[$b)$] a local minimizer $\bfal^*$ is an $\eta$-stationary point either for any $\eta>0$ if  $\|\bfal^*\|_0<s$ or for any $0<\eta\leq\eta^*$   if  $\|\bfal^*\|_0=s$, where $\eta^*>0$ is relied on $\bfal^*$.
\end{itemize}
ii)  For global  minimizers and $\eta$-stationary points, it has
\begin{itemize}
\item[$c)$] if $\|\bfal^*\|_0<s$, then the local minimizer, the  global minimizer and the $\eta$-stationary point are identical to each other and unique.
\item[$d)$] if $\|\bfal^*\|_0=s$, and the point $\bfal^*$ is an $\eta$-stationary point  with  $\eta\geq C$, 
then it is also a global minimizer. Moreover, it is also unique if $\eta> C$.
\end{itemize} 
\end{theorem}
 {\cref{existence-global-minimizers} states that the global minimizers (also being a local minimizer) of \eqref{SM-SVM-h-equ-dual-sparse} exist. In addition, \cref{nec-suff-opt-con-tau} b) shows that a local minimizer of \eqref{SM-SVM-h-equ-dual-sparse} is an $\eta$-stationary point.   Therefore, the $\eta$-stationary point exists,  indicating that there must exist $b^*$ such \eqref{SM-SVM-h-equ-dual-sparse-tau}.}
\section{Newton Method}\label{sec:SNM}
This section applies the Newton method for solving  equation \eqref{stationary-equation-T}.  Let $\bfz^k$ be defined in \eqref{z-z*-zk} and the current approximation to a solution of (\ref{stationary-equation-T}). Choose one index set 
\begin{equation}\label{Newton-Method-Tk} T_k\in\T_s(\bfal^k-\eta g(\bfz^k)).\end{equation}
Then find the Newton direction  $\bfd^{k}\in\R^{m+1}$ by solving the following linear equations,
\begin{equation} \label{Newton-Method}
\nabla F (\bfz^k;T_k) \bfd^k  = - F(\bfz^k;T_k).
\end{equation}
Substituting \eqref{stationary-equation-T} and \eqref{Jacobian-F} into \eqref{Newton-Method}, we   obtain
\begin{eqnarray}\label{sequation-k-0}
\left[
\begin{array}{ccc}
 H_{T_k}(\bfal^k) &0 & \bfy_{T_k}\\ 
  0&I&0 \\
 \bfy_{T_k}^\top&0&0
\end{array}\right]\left[
\begin{array}{l}
 \bfd^k_{T_{k}} \\ 
 \bfd^{k}_{ \overline{ T }_{k}} \\
 d^k_{m+1}
\end{array}\right] = -\left[
\begin{array}{c}
 g_{T_k} (\bfz^k)\\ 
  \bfal^{k}_{\overline{ T }_{k}} \\
\langle\bfal_{T_k}^k, \bfy_{T_k} \rangle
\end{array}\right].
\end{eqnarray}
After we get the direction, the full Newton step size is taken, namely, $\bfz^{k+1}=\bfz^k+\bfd^k$, to derive
\begin{eqnarray}
\label{form4}
\bfz^{k+1}=\left[
\begin{array}{l}
 \bfal^{k}_{T_{k}} \\ 
  \bfal^{k}_{ \overline{ T }_{k}} \\
   b ^{k}
\end{array}\right]+\left[
\begin{array}{l}
 \bfd^k_{T_{k}} \\
 \bfd^{k}_{ \overline{ T }_{k}} \\
 d^k_{m+1}
\end{array}\right] 
=\left[
\begin{array}{c}
 \bfal^{k}_{T_{k}} + \bfd^k_{T_{k}} \\ 
  0 \\
   b ^{k}+d^k_{m+1}
\end{array}\right]. 
\end{eqnarray} 
 Now we summarize the whole framework in  \cref{alg:snsvm}. 
\begin{algorithm}[H]
\caption{Newton method for sparse SVM.}
\label{alg:snsvm}
\begin{algorithmic}
\STATE{Give parameters $C>c,\eta, \epsilon, K>0, s\in[m]$. }
\STATE{Initialize $\bfz^0$, pick $T_{0}\in\T_s(\bfal^0-\eta g(\bfz^0))$ and set $k:=0$.}
\WHILE{$\| F (\bfz^k;T_k)  \| \geq \epsilon$ and $k\leq K$}
\STATE{ Update $\bfd^k$ by solving (\ref{sequation-k-0}).}
\STATE{ Update  $\bfz^{k+1}$ by (\ref{form4}).}
\STATE{ Update $T_{k+1}\in\T_s(\bfal^{k+1}-\eta g(\bfz^{k+1}) )$ and  $k:=k+1$.}
\ENDWHILE
\RETURN the solution $\bfz^{k}$.
\end{algorithmic}
\end{algorithm}
It is easy to see from \eqref{sequation-k-0} that $\bfd^{k}_{ \overline{ T }_{k}}=-\bfal^{k}_{\overline{ T }_{k}}$ is derived directly. Therefore, every new point is always sparse due to $$\|\bfal^{k+1}\|_0\overset{\eqref{form4}}{=}\|\bfal^{k}_{T_{k}} + \bfd^k_{T_{k}}\|_0 \leq  |T_k|=s.$$ Moreover, the major computation in \eqref{sequation-k-0} is from the part on $T_k$.  However,  $|T_k|=s$ can be controlled to have a very small scale compared to $m+1$,  leading to a considerably low computational complexity.
 \subsection{Complexity  Analysis}
To derive the Newton direction in \eqref{sequation-k-0}, we need to calculate 
\begin{eqnarray}\label{sequation-k-euq}
d^k_{m+1}  &=& -\frac{\left\langle \bfy_{T_k},~\Theta^{-1}  g_{T_k}(\bfz^{k})-\bfal_{T_k}^k  \right\rangle }{ \langle \bfy_{T_k},~\Theta^{-1}\bfy_{T_k}^k  \rangle} ,\nonumber\\
 \bfd^{k}_{   { T }_{k}~~}  &= & - \Theta^{-1}  \left[g_{T_k}(\bfz^{k})+ d^k_{m+1}\bfy_{T_k} \right],\\ 
 \bfd^{k}_{ \overline{ T }_{k}}~~ &=&-\bfal^{k}_{\overline{ T }_{k}},\nonumber
\end{eqnarray}
where $\Theta:=H_{T_k}(\bfal^k)$. For the computational complexity of  \cref{alg:snsvm},  we  observe that calculations of $\Theta^{-1}$ and $\T_s$ dominate the whole computation. Their  total complexity in each step is approximately
\begin{eqnarray}
\label{complexity-total} {O}\left(mn+\max\{n,s\}s^2\right).\end{eqnarray} 
This can be derived by the following analysis: 
\begin{itemize}
\item Recall the definition \eqref{gradient-hessian} of $H(\cdot)$ that 
$$\Theta=(E(\bfal^k))_{T_kT_k}+Q_{T_k}^\top Q_{T_k}.$$
The complexity of computing $\Theta$ is about ${O}(ns^2)$ since $|T_k|=s$. Computing its inverse has the complexity of at most ${O}(s^3)$.  Therefore, the complexity of deriving $\Theta^{-1}$ is $ {O}(\max\{n,s\}s^2).$
\item  To pick $T_{k+1}$ from $\T_s$, we need to compute $g(\bfz^{k+1})$
 and  select the $k$ largest elements of $|\bfal^{k+1}-\eta g(\bfz^{k+1})|$. The complexity of computing the former is $O(mn)$ due to the computation of  $Q^\top (Q_{T_k}\bfal^{k+1}_{T_k})$ and the complexity of computing the latter is $O(m + s {\rm ln} s)$. Here, we use  the  MATLAB built-in function {\tt maxk}  to select the $s$ largest elements.  
\end{itemize} 
%According to  \eqref{complexity-total}, if $s$ is chosen to be quite small, (e.g., $s = 0.01m$), the complexity in \eqref{complexity-total} can be significantly reduced compared to the  complexity of solving the kernel-based SVM without applying data reduction strategies. Hence, this enables the extremely large scale computation.
 {We now summarize the complexity of \cref{alg:snsvm} and some other methods in  \cref{table:com}, where in  \cite{LM01} and \cite{TLT10,LEP19}, their proposed method solved  a quadratic programming  in each iteration $k$, and $N_k$ was the number of iterations used to solve the quadratic programming. In \cite{FS01, ZSW20}, $r$ was the rank of $Q^\top Q$ or its an approximation.  Moreover, $\bar m$ was the number of the pre-selected samples from the total samples \cite{LM01}, and $m_k$ was the number of the samples used at the  $k$th iteration \cite{NMTH10}. In the setting where data involves large numbers of samples, i.e., $m$ is considerably large in comparison with $n$, the rank $r\approx n$, which means \cref{alg:snsvm} has the lowest computational complexity among the second-order methods. It is well-known that second-order methods can converge quickly within a small number of iterations. Therefore, \cref{alg:snsvm} can be executed super-fast, which has been demonstrated by numerical experiments.
%Moreover,    the first three methods have a similar computational complexity that is approximately $\mathcal{O}(mn)$.
\begin{table}[H] 
\centering
\caption{Complexity of different algorithms.}
\label{table:com} 
\renewcommand{\arraystretch}{1.25}\addtolength{\tabcolsep}{-2pt}	
\begin{tabular}{p{1.5cm}llr}\hline
Methods	&Reference&  Complexity &	 Descriptions \\\hline
\multirow{5}{*}{First-order}
&\cite{LM01} & $\mathcal{O}(mn+N_k\bar m^2)$& $\bar m\in[1,m]$ \\
&\cite{CSS13} & $\mathcal{O}(mn)$&  \\
&\cite{NMTH10} & $\mathcal{O}(m_k^3+m m_k)$&  $m_k\in[1,m]$\\
&\cite{TLT10} &  $\mathcal{O}(N_kmn+n{\rm ln}(B))$& $B\in(1,n]$ \\&\cite{LEP19}  & $\mathcal{O}(N_k m^2)$&\\\hline
\multirow{4}{*}{Second-order}&
\cite{FS01} &$\mathcal{O}(m r^2)$&$r\in[1,m]$\\
&\cite{ZSW20}&$\mathcal{O}(m^2(n+r))$&$r\in[1,m]$\\
&\cite{LDD11} & $\mathcal{O}(m^3)$&  \\
&\cref{alg:snsvm} & ${O}\left(mn+\max\{n,s\}s^2\right)$&$s \in[1,m]$\\
\hline
\end{tabular}
\end{table}}
%For the method   RSVM   in \cite{LM01}, the computation comprises of calculating a matrix production with the complexity  $\mathcal{O}(mn)$ and addressing a quadratic programming with $\bar m\in[1,m]$ variables. The total one is  $\mathcal{O}(mn+\bar m^2)$. In \cite{CSS13}, each step updating the subgradient and needs the complexity $\mathcal{O}(mn)$. As stated in \cite{TLT10}, their method FGM has the complexity  $\mathcal{O}(mn+n{\rm ln}(B))$, where $B\in(1,n]$. In \cite{NMTH10}, the condensed AVM training algorithm takes the complexity  $\mathcal{O}(m_k^3+m m_k)$, where $m_k$ is the number of support vectors used in $k$th iteration. The algorithm ISLS-SVM proposed in \cite{LDD11} solves an linear equation with size $(m+1)$ and calculate a product between two kernel matrices matrix with the order $m\times m$ in each step. Therefore, the complexity is almost $\mathcal{O}(m^3)$. The $L_0$ SVM Algorithm developed in \cite{LEP19} needs to address an quadratic programming with variable size $m$ in each iteration, the complexity is then $\mathcal{O}(m^2)$. }
 
 \subsection{One-Step Convergence}
 For a point  $
\bfz^*=(\bfal^*; b ^*)$ with $\bfal^*$ feasible to  \eqref{SM-SVM-h-equ-dual-sparse}, define
 \allowdisplaybreaks\begin{eqnarray}
\label{tau*}\eta^*&:=&\begin{cases}
  \|\bfal^*\|_{[s]} \|g(\bfz^*)\|_{[1]}^{-1} , &\| \bfal^*\|_0=s,\\
+\infty, & \| \bfal^*\|_0<s,
\end{cases}
\end{eqnarray}
The convergence result is stated by the following theorem. 
\begin{theorem}\label{the:i-iterate-convergence}
Let $\bfz^*$ be an $\eta$-stationary point of the problem \eqref{SM-SVM-h-equ-dual-sparse} with $0<\eta<\eta^*$, where $\eta^*$ is given by \eqref{tau*}. Let $\{\bfz^k\}$ be the sequence generated by  \cref{alg:snsvm}. There always exists a $\delta^*>0$ such that if at a certain iteration,  $\bfz^{k}\in N(\bfz^*,\delta^*)$,  then 
 $$\bfz^{k+1}=\bfz^*,~~~~\|F(\bfz^{k+1},T_{k+1})\|=0.$$ Namely, \cref{alg:snsvm} terminates at the th$(k+1)$ step. 
\end{theorem}
 \begin{figure}[!th]  
\centering
\begin{subfigure}{.24\textwidth}
	\centering
	\includegraphics[width=1.04\linewidth]{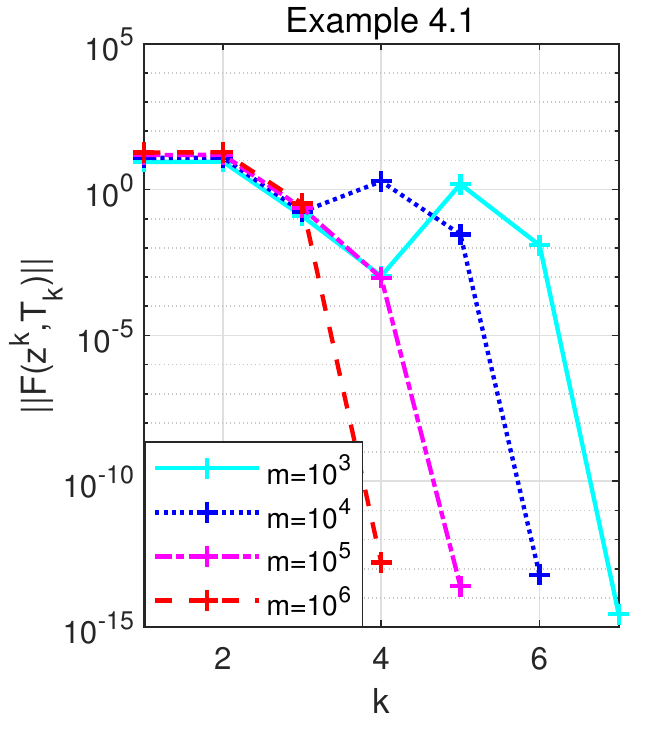}
	%\caption{$m=200$}
	%\label{fig:USairport2010_result2}
\end{subfigure}
\begin{subfigure}{.24\textwidth}
	\centering
	\includegraphics[width=1.04\linewidth]{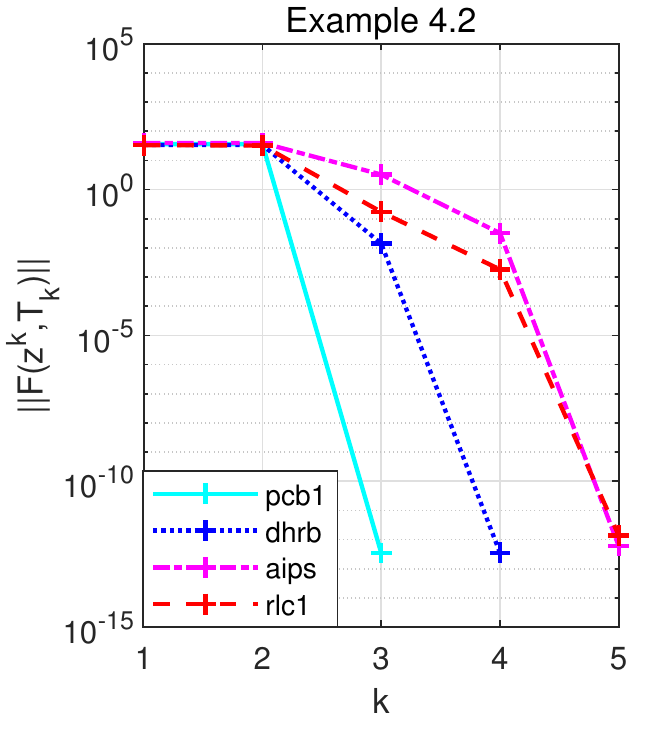}
	%\caption{$m=400$}
	%  \label{fig:Blogs_result2}
\end{subfigure} \\\vspace{-1mm}
\caption{Decreasing of $\|F(\bfz^{k},T_{k})\|$.}
\label{fig:F-iter}
\end{figure}

We note that  \cref{alg:snsvm} can terminate at the next step if the current point falls into a local area of an $\eta$-stationary point. This means that if by chance the starting point is chosen within the local area,  the proposed algorithm will take one step to terminate. Hence, it enjoys a fast convergence property.  {Of course, this depends on the choices of the initial point. In reality, how close the initial point is to an $\eta$-stationary point is unknown. Nevertheless,  numerical experiments show that the proposed method is insensitive the choice of the starting point. It can take a few steps to converge. For instance,  the results of \cref{alg:snsvm} solving \cref{ex:syn-data} and \cref{ex:real-data} with the initial point $\bfal^0=0$ and $ b ^0={\rm sgn}(\langle\bfy, \bfe\rangle)$ were presented in Fig. \ref{fig:F-iter}. The line of {\tt pcb1} showed that the point $\bfz^{3}$ at step 3 may fall into a local area of an $\eta$-stationary point since $\|F(\bfz^{4},T_{4})\|(\approx 10^{-12})$ almost reached  0 at step 4. Therefore, it terminated at step 4. This showed potential of termination at the next step, even though the initial point was not close to the $\eta$-stationary point. Similar observations can be found for all other lines.}

 {Moreover, \cref{nec-suff-opt-con-tau} states that the $\eta$-stationary point is a local minimizer of the problem \eqref{SM-SVM-h-equ-dual-sparse}. Hence, \cref{alg:snsvm} at least converges to a local minimizer. If further $\|\bfz^*\|_0 < s$ or $\|\bfz^*\|_0 = s$ and $\eta^*\geq C$, then the $\eta$-stationary point is also a global minimizer by \cref{nec-suff-opt-con-tau} c) and d), namely, \cref{alg:snsvm}  converges to the global minimizer. }

\subsection{Sparsity Level Tuning} 
One major issue encountered in practice is that the sparsity level $s$ in \eqref{SM-SVM-h-equ-dual-sparse} is usually unknown beforehand. The level has two important effects: (i)  A larger $s$ corresponds to better classifications since more samples are taken into consideration. (ii) However, a smaller $s$ leads to higher computational speed of  \cref{alg:snsvm} because the complexity in \eqref{complexity-total} depends on $s$. Moreover, according to \eqref{Representer-Theorem} that the number of support vectors is smaller than $\|\bfal\|_0\leq s$, a smaller $s$ results in a smaller number of support vectors.   Therefore, to balance these two requirements, we design the following rule to update the unknown $s$. We start with a small integer and then increase it until  the following halting conditions are satisfied
\begin{eqnarray}  
\label{halt-cods}   \| F (\bfz^k;T_k)  \|&<& \epsilon,\\
\label{halt-cods-1} \Big|{\tt ACC}(\bfal^k)-\max_{j\in[ k-1]}{\tt ACC}(\bfal^j)\Big|&<&10^{-4},\end{eqnarray}
where we admit ${\tt ACC}(\bfal^{-1})=0 $ and  ${\tt ACC}(\bfal)$ is defined by
\begin{eqnarray}\label{AC}
{\tt ACC}(\bfal):=\left[1-  \frac{1}{m}{\|{\rm sgn}(X\bfal+ b )-\bfy\|_0}\right]\times 100\%, 
\end{eqnarray}
with ${\rm sgn}(t) = 1$ if $t>0$ and $-1$ otherwise. The halting condition \eqref{halt-cods} means that $\bfal^k$ is almost an $\eta$-stationary point due to the small  $\| F (\bfz^k;T_k)  \|$, while \eqref{halt-cods-1} implies that the classification accuracy ${\tt ACC}(\bfal^k)$ does not increase significantly. Therefore, it is unreasonable to further increase $s$ to achieve a better solution because a larger $s$ will lead to greater computational cost. %In other words, it is reasonable to terminate the algorithm if  $\bfal^k$ satisfies such conditions. 
Our numerical experiments demonstrate that \cref{alg:snsvm} under this scheme works very well.  Overall, we derive {\tt NSSVM} in \cref{alg:NSSVM}. 
\begin{algorithm}[!th]
\caption{{\tt NSSVM}: Newton method for sparse SVM with adaptively tuning the sparsity level $s$.}
\label{alg:NSSVM}
\begin{algorithmic}
\STATE{Give parameters $C>c,\eta, \epsilon, K>0, \sigma>1, s_0\in[m]$}. 
\STATE{Initialize $\bfz^0$, pick $T_{0}\in\T_{s_0}(\bfal^0-\eta g(\bfz^0))$ and set $k:=0$.}
\WHILE{$\bfz^k$ violates \eqref{halt-cods} and \eqref{halt-cods-1} and $k\leq${\tt MaxIt}}
\STATE{ Update $\bfd^k$ by solving (\ref{sequation-k-0}).}
\STATE{ Update  $\bfz^{k+1}$ by (\ref{form4}).}
%\STATE{ Update ${\tt MaxACC}=\max\{{\tt ACC}(\bfal^1),\cdots,{\tt ACC}(\bfal^{k-1})\}.$ }
\STATE{ Update $s_{k+1}= 
\sigma s_k$ if $k$ is a multiple of $10$ or \eqref{halt-cods} is met, and $s_{k+1}= s_k$ otherwise.}
%$s_{k+1}=\begin{cases} 
%\sigma s_k,& \text{if \eqref{halt-cods} is met or $k$ is a multiple of $10$},\\
%s_k,& \text{otherwise}.
%\end{cases}$}
\STATE{ Update $T_{k+1}\in\T_{s_{k+1}}(\bfal^{k+1}-\eta g(\bfz^{k+1}) )$ and $k:=k+1$.}
\ENDWHILE
\RETURN the solution $\bfz^{k}$.
\end{algorithmic}
\end{algorithm}

\section{Numerical Experiments}\label{numerical}

This part conducts extensive numerical experiments of  \cref{alg:snsvm} and  \cref{alg:NSSVM} ({\tt NSSVM}\footnote{\url{https://github.com/ShenglongZhou/NSSVM}}) by using MATLAB (R2019a) on a laptop with  $32$GB memory and Inter(R) Core(TM) i9-9880H 2.3Ghz CPU.
\subsection{Testing Examples}
We first consider a two-dimensional example with synthetic data, where the features come from Gaussian distributions.
\begin{example}[Synthetic data in $\R^2$ \cite{XAC2016,  HSS16}]\label{ex:syn-data} Samples $\bfx_{i}$ with positive labels $y_i=+1$ are drawn from  the normal distribution with mean $(0.5,-3)^{\top}$ and variance $\Sigma$, and samples  $\bfx_{j}$ with negative labels $y_j=-1$ are drawn from the normal distribution with mean $(-0.5,3)^{\top}$ and variance $\Lambda$, where $\Sigma$ and $\Lambda$ are diagonal matrices with $\Sigma_{11}=\Lambda_{11}=0.2,~\Sigma_{22}=\Lambda_{22}=3.$  We generate $2m$ samples with equal numbers of two classes and then evenly split them into a training and a testing set. %Finally, we randomly flip $rm$ labels in the training data, namely, $rm$ samples are treated as outliers.
\end{example}

\begin{example}[Real data in higher dimensions]\label{ex:real-data} We select 30 datasets with $m \gg n$  from the libraries: libsvm\footnote{\url{https://www.csie.ntu.edu.tw/~cjlin/libsvmtools/datasets/}}, uci\footnote{\url{http://archive.ics.uci.edu/ml/datasets.php}} and kiggle\footnote{\url{https://www.kaggle.com/datasets}}. Here, $aeb:=a\times10^b$. All datasets are feature-wisely scaled to $[-1,1]$ and all of the classes that are not $1$ are treated as $-1$.  Their details are presented in  \cref{Table-svm-more-m-less-n}, where $7$ datasets have the testing data. For each of the dataset without the testing data, we split the dataset into two parts. The first part contains 90\% of samples treated as the training data and the rest are the testing data.
\end{example}

\begin{table}[!th] 
	\renewcommand{\arraystretch}{.95}\addtolength{\tabcolsep}{-3.5pt}
	\caption{Descriptions of real datasets.}\vspace{-6mm}
	\label{Table-svm-more-m-less-n}
	\begin{center}
		\begin{tabular}{lllrrrr }  \\	\hline
	Data	&	Datasets		&	Source	&	 		&	Train	&	Test\\ 
		 	&	 	&	 	&	$n$		&	$m$	&	$m_t$\\ \hline
		 	 \texttt{pcb1}&	\multirow{5}{3cm}{Polish companies bankruptcy data}	&	uci	&	64 	&	7027 	&	0	&	\\
 \texttt{pcb2}&		&	uci	&	64 	&	10173 	&	0	&	\\
 \texttt{pcb3}&		&	uci	&	64 	&	10503 	&	0	&	\\
 \texttt{pcb4}&		&	uci	&	64 	&	9792 	&	0	&	\\
 \texttt{pcb5}&		&	uci	&	64 	&	5910 	&	0	&	\\
  \texttt{a5a}&	A5a	&	libsvm	&	123 	&	6414 	&	26147	&		\\
 \texttt{a6a}&	A6a	&	libsvm	&	123 	&	11220 	&	21341	&		\\
 \texttt{a7a}&	A7a	&	libsvm	&	123 	&	16100 	&	16461	&		\\
 \texttt{a8a}&	A8a	&	libsvm	&	123 	&	22696 	&	9865	&		\\
 \texttt{a9a}&	A9a	&	libsvm	&	123 	&	32561 	&	16281	&	\\\hline
 \texttt{mrpe}&	malware analysis datasets	&	kaggle	&	1024 	&	51959 	&	0	\\
  \texttt{dhrb}&	hospital readmissions bnary	&	kaggle	&	17 	&	59557 	&	0	\\
\texttt{aips}&	airline passenger satisfaction	&	kaggle	&	22 	&	103904 	&	25976	\\
 \texttt{sctp}&	santander customer transaction	&	kaggle	&	200 	&	200000 	&	0	\\
 \texttt{skin}&	skin\_nonskin	&	libsvm	&	3 	&	245056 	&	0	\\
 \texttt{ccfd}&	credit card fraud dtection	&	kaggle	&	28 	&	284807 	&	0	\\
 \texttt{rlc1}&	 \multirow{3}{4cm}{record linkage comparison patterns}	&	uci	&	9 	&	574914 	&	0	\\
% \texttt{rlc2}&		&	uci	&	9 	&	574914 	&	0	\\
 $\cdots$&		&	 $\cdots$	&	 $\cdots$ 	&	 $\cdots$ 	&	 $\cdots$	\\
% \texttt{rlc4}&		&	uci	&	9 	&	574914 	&	0	\\
% \texttt{rlc5}&		&	uci	&	9 	&	574914 	&	0	\\
% \texttt{rlc6}&		&	uci	&	9 	&	574914 	&	0	\\
% \texttt{rlc7}&		&	uci	&	9 	&	574914 	&	0	\\
% \texttt{rlc8}&		&	uci	&	9 	&	574914 	&	0	\\
% \texttt{rlc9}&		&	uci	&	9 	&	574914 	&	0	\\
 \texttt{rlc10}&		&	uci	&	9 	&	574914 	&	0	\\ 
 \texttt{covt}&	covtype.binary	&	libsvm	&	54 	&	581012 	&	0	\\
% \texttt{retb}&	real time bidding	&	kaggle	&	88 	&	$1e6$ 	&	0	\\
 \texttt{susy}&	susy	&	uci	&	18 	&	$5e6$ 	&	0	\\
 \texttt{hepm}&	hepmass	&	uci	&	28 	&	$7e6$ 	&	$35e5$	\\
 \texttt{higg}&	higgs	&	uci	&	28 	&	$11e6$ 	&	0	\\
 \hline
		\end{tabular}
	\end{center}
\end{table}
To compare the performance of all selected methods, let $\bfal$ be the solution/classifier generated by one method. We report  the  CPU time ({\tt TIME}), the training classification accuracy ({\tt ACC}) by
\eqref{AC} where $X$ and $\bfy$ are the training samples and classes, the testing classification accuracy ({\tt TACC}) by
\eqref{AC} where $X$ and $\bfy$ are the testing samples and classes, and the number of support vectors ({\tt NSV}).
\subsection{Implementation and parameters tuning}
The starting point $
\bfz^0$ is initialized as $\bfal^0=0$ and $ b ^0={\rm sgn}(\langle\bfy, \bfe\rangle)$ if no additional information is provided.   The maximum number of iterations and the tolerance are set as $K=1000$, $\epsilon =\max\{\sqrt{m}$ and $\sqrt{n}\}10^{-6}$. Recall that the model \eqref{SM-SVM-h-equ-dual-sparse} involves three important parameters  $C, c$ and $s$ and the $\eta$-stationary point has the parameter $\eta$. We now apply \cref{alg:snsvm} for tuning them as described below.

\subsubsection{Effect of $\eta$}
\begin{figure}[!th]  
\centering
\begin{subfigure}{.24\textwidth}
	\centering
	\includegraphics[width=1.04\linewidth]{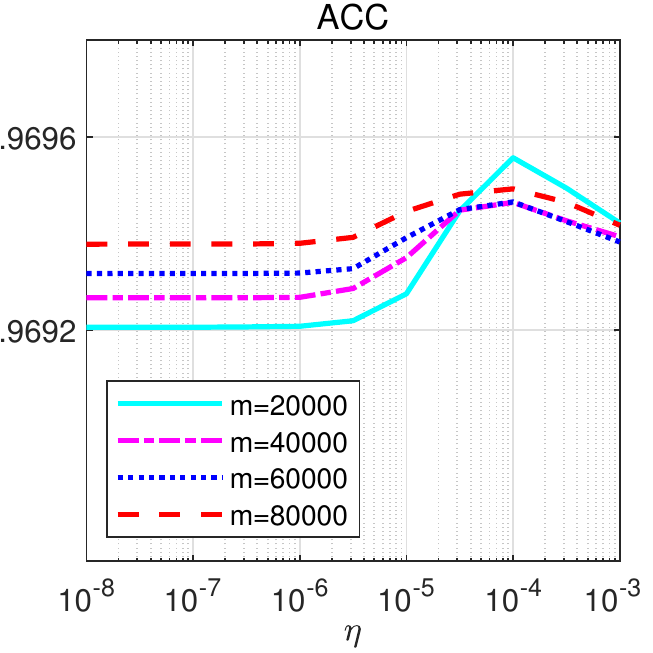}
	%\caption{$m=200$}
	%\label{fig:USairport2010_result2}
\end{subfigure}
\begin{subfigure}{.24\textwidth}
	\centering
	\includegraphics[width=1.04\linewidth]{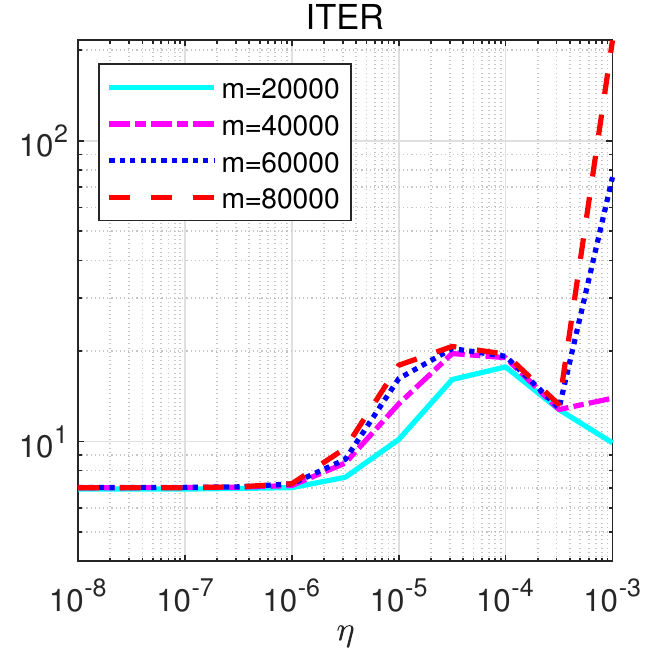}
	%\caption{$m=400$}
	%  \label{fig:Blogs_result2}
\end{subfigure} \\\vspace{-1mm}
\caption{Effect of  $\eta$.}
\label{fig:effect-eat}
\end{figure}

For \cref{ex:syn-data}, we fix  $C=2^{-2}, c=C/2, s=0.005m$ and vary $\eta\in[10^{-8},10^{-2}]$ and $m\in\{2,4,6,8\}\times 10^4$ to examine the effect of $\eta$ on \cref{alg:snsvm}. For each case $(m,\eta)$, we run 100 trials  and report the average results in Fig. \ref{fig:effect-eat}. It is clearly observed that  the  method was insensitive to $\eta$   in the range $[10^{-8},10^{-6}]$. For each $m$, the  {\tt ACC} lines  were stabilized at a certain level when $\eta\in[10^{-8},10^{-6}]$,   increased when $\eta\in[10^{-6},10^{-4}]$ and then decreased. It appears that $\eta$ of approximately $10^{-4}$ gave the best performance for  \cref{alg:snsvm} in terms of the highest {\tt ACC}.
%We observed that some $\eta$ (e.g., $10^{-3}$) led  to a large number of iterations for the algorithm. This was because the large  $\eta$ greatly varies the $T_k$ determined by \eqref {Newton-Method-Tk}, leading to the slow convergence of the algorithm. 
Moreover, we tested problems with different dimensions and found that $\eta=1/m$ enabled the algorithm to achieve steady behaviour. Therefore, we set this option for $\eta$ in the subsequent numerical experiments if no additional information is provided.

\subsubsection{Effect of  $C$ and $c$}
\begin{figure}[H]  
\centering
\begin{subfigure}{.24\textwidth}
	\centering
	\includegraphics[width=1.05\linewidth]{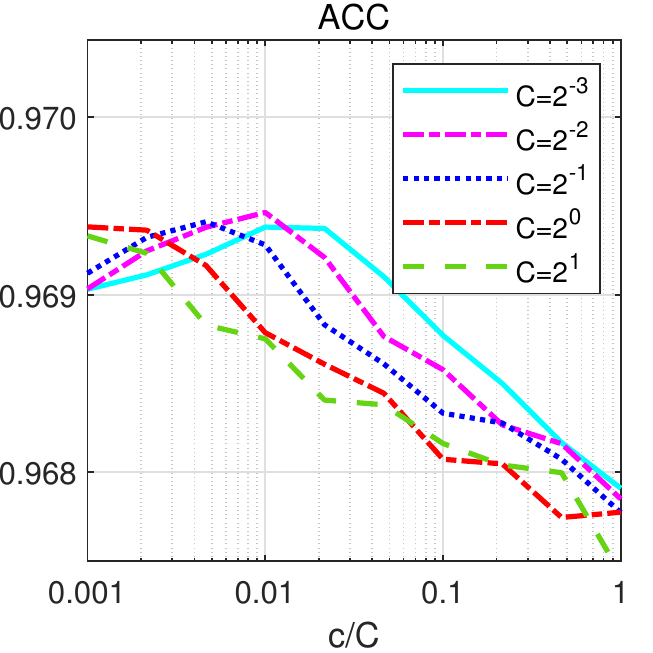}
	%\caption{$m=200$}
	%\label{fig:USairport2010_result2}
\end{subfigure}
\begin{subfigure}{.24\textwidth}
	\centering
	\includegraphics[width=1.05\linewidth]{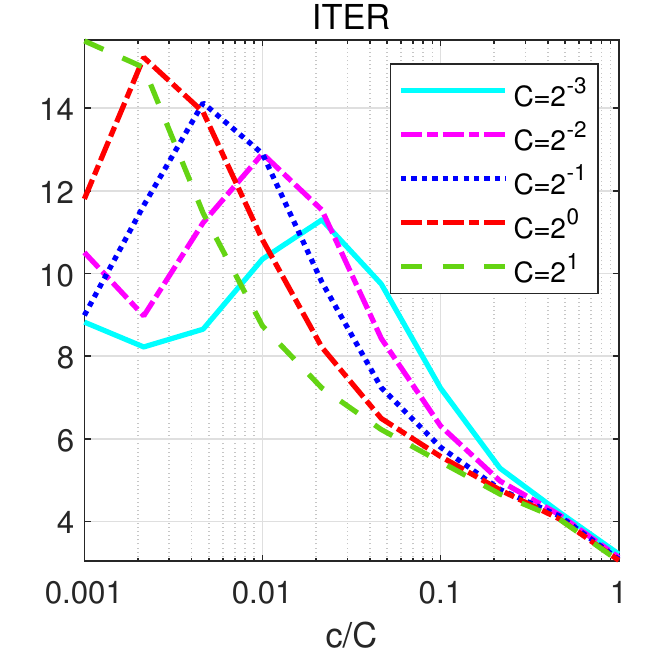}
	%  \label{fig:Blogs_result2}
\end{subfigure} \\\vspace{-1mm}
\caption{Effect of  $C$ and $c$.}
\label{fig:effect-cc}
\end{figure}

To see the effect of  parameters $C$ and $c$, we choose $C\in\{2^{-3},2^{-2},2^{-1},2^{0},2^{1}\}$ and $c=a C$, where $a\in[0.001,1]$. Again, we apply \cref{alg:snsvm}  for solving   \cref{ex:syn-data} with fixing $m=10000$ and $s=0.005m$. The average results were recoded in Fig. \ref{fig:effect-cc}, where it was observed that {\tt ACC} reached the highest peak at $a=c/C=0.01$. When $c/C\geq0.02$, the bigger $C$ was, the higher {\tt ACC} was but the larger {\tt ITER} was as well.  
To balance the accuracy and the number of iterations, in the following numerical experiments, we set  $C=0.25$ and $c=0.01C$   unless specified otherwise.

\subsubsection{Effect of $s$}
\begin{figure}[H]  
\centering
\begin{subfigure}{.24\textwidth}
	\centering
	\includegraphics[width=1.05\linewidth]{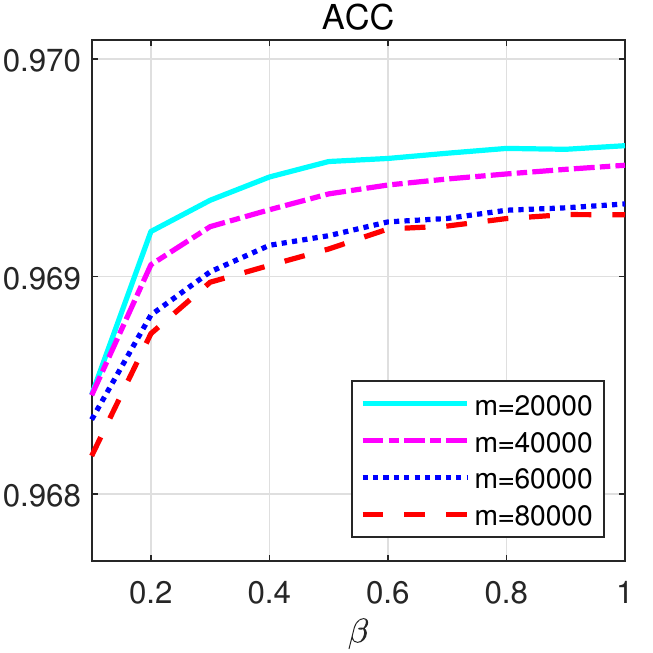}
	%\caption{$m=200$}
	%\label{fig:USairport2010_result2}
\end{subfigure}
\begin{subfigure}{.24\textwidth}
	\centering
	\includegraphics[width=1.05\linewidth]{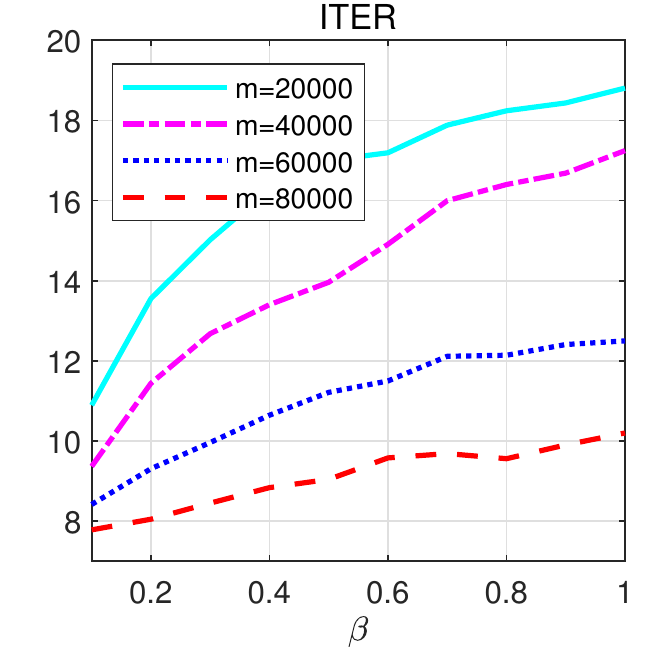}
	%  \label{fig:Blogs_result2}
\end{subfigure} 
%\\\vspace{1mm}
%\begin{subfigure}{.24\textwidth}
%	\centering
%	\includegraphics[width=1.05\linewidth]{TIME-s-eps-converted-to.pdf}
%	%  \label{fig:Blogs_result2}
%\end{subfigure}
%\begin{subfigure}{.24\textwidth}
%	\centering
%	\includegraphics[width=1.05\linewidth]{NSV-s-eps-converted-to.pdf}
%	%  \label{fig:Blogs_result2}
%\end{subfigure} \vspace{-1mm}
\caption{Effect of  $s$.}
\label{fig:effect-s}
\end{figure}

Apparently, the sparsity level $s$ has a strong influence on the  solutions to the problem \eqref{SM-SVM-h-equ-dual-sparse} and impacts the computational complexity of  \cref{alg:snsvm}.  Below, we use the following rule to tune the sparsity level $s$,
$$s(\beta):=\lceil\beta n({\rm log}_2(m/n))^2\rceil,$$
where $ \lceil t\rceil$ presents the smallest integer that is no less than $t$, and $\beta$ is chosen based on the problem to be solved. We use this rule because it relies on the dimensions of a given problem and allows \cref{alg:snsvm} to deliver an overall desirable performance. By setting $s=s(\beta)$, we focus on selecting a proper $\beta$. We note that $s$ increases with the rising of $\beta$. As demonstrated in Fig. \ref{fig:effect-s}, where  \cref{ex:syn-data} was solved again, the average results showed the effect of $\beta\in[0.1,1]$ (i.e., the effect of $s$) of the algorithm. %Based on the model \eqref{SM-SVM-h-equ-dual-sparse}, $s$ is  the number of support vectors, namely, ${\tt NSV} =s$, see the last plot in \cref{fig:effect-s}.  
It is observed from the figure that  {\tt ACC} increased with increasing $s$ but gradually stabilized at a certain level after $\beta>0.4$, indicating that there was no need to increase $s$ beyond this point.  Furthermore, the small number of iterations showed that \cref{alg:snsvm} can converge quickly.

% It is observed from the figure that  {\tt ACC} increased with increasing $s$ but gradually stabilized at a certain level after $\beta>0.4$, indicating that there was no need to increase $s$ beyond this point.  Furthermore, the small number of iterations indicated \cref{alg:snsvm} can converge quickly.
%Furthermore, the number of iterations and computational time ascended with the rising of $s$. Finally,  \cref{alg:snsvm} ran quite fast, solving each instance in less than 0.02 seconds. This high speed was achieved due to the low computational complexity of \eqref{complexity-total} in each iteration and a small number of iterations required to achieve convergence (e.g., less than 20 for all cases).
\subsubsection{Effect of the initial points}
\begin{figure}[!th]  
	\centering
	\includegraphics[width=1.01\linewidth]{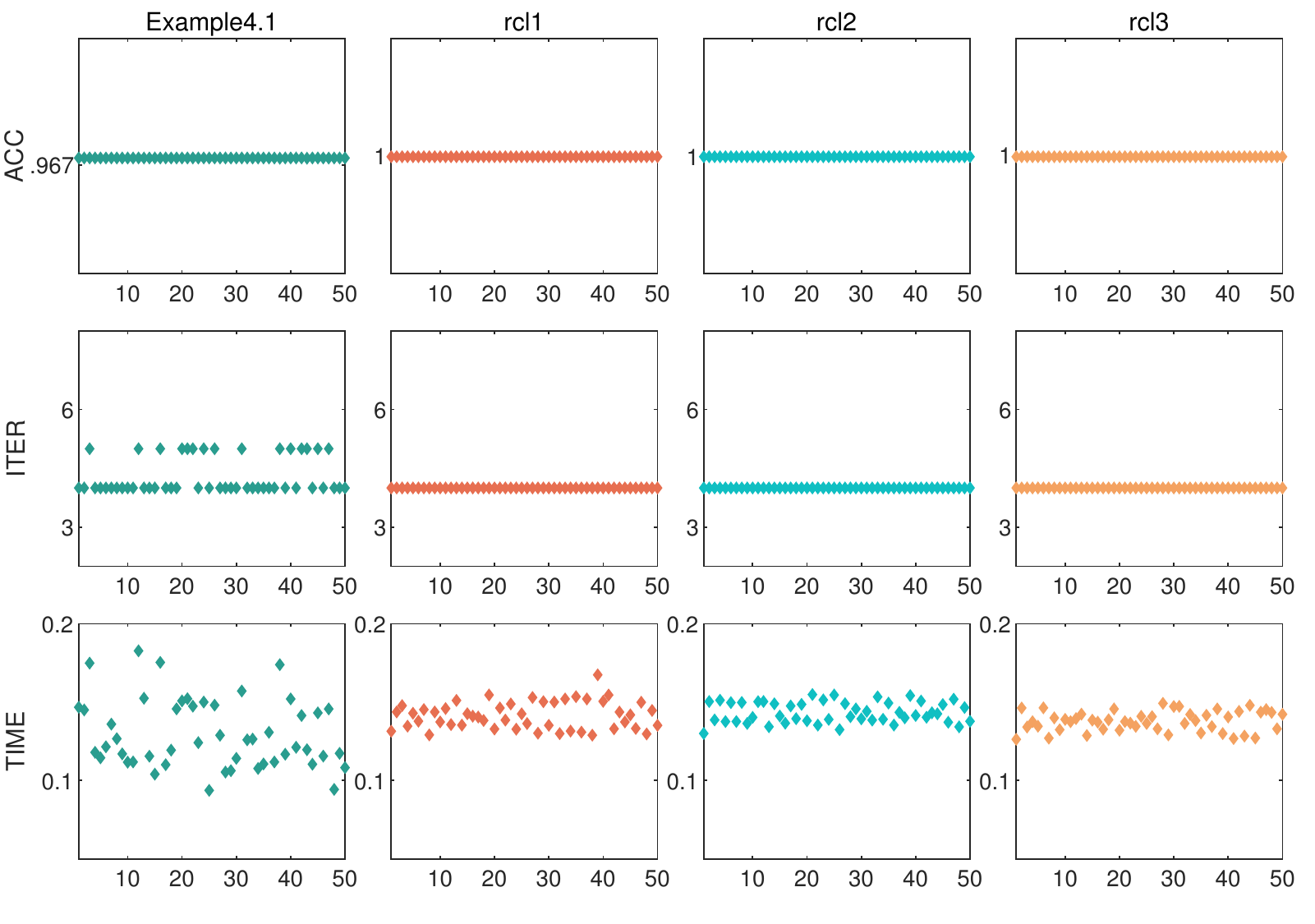}
	%\caption{$m=200$}
	%\label{fig:USairport2010_result2}
\caption{Robust to  the initial points $\bfal^0$.}
\label{fig:effect-alpha0}
\end{figure}

 \cref{the:i-iterate-convergence} states that the algorithm can converge at the next step if the initial point is sufficiently close to an $\eta$-stationary point that however is unknown beforehand. Therefore, we examine how the initial points affect the performance of \cref{alg:snsvm}. To carry out this study, we apply the algorithm for  solving \cref{ex:syn-data} with $m=10^6$ and $s=s(1)$, and \cref{ex:real-data} with datasets: {\tt rcl1}, {\tt rcl2}, {\tt rcl3} and $s=s(0.05)$. For each data, we run the algorithm under $50$ different initial points $\bfal^0$. The first point is $\bfal^0=0$ and the other $49$ points have randomly generated  entries from the uniform distribution, namely, $\alpha_i^0\sim U[0,1]$. The results were plotted in Fig. \ref{fig:effect-alpha0}, where the x-axis represents the 50 initial points. For each dataset,  all  {\tt ACC} stabilized at a certain level, indicating that these initial points  had little effect on the accuracy. The results for ${\tt ITER}$ illustrated that the algorithm converged very quickly, stopping within $6$ steps for all cases. To summarize,  \cref{alg:snsvm} was not too sensitive to the choice of the initial point for these datasets. Hence, for simplicity, we initialize our algorithm with $\bfal^0=0$ and $ b ^0={\rm sgn}(\langle\bfy, \bfe\rangle)$  
 in the sequel.

\subsection{\cref{alg:snsvm} v.s. \cref{alg:NSSVM}: {\tt NSSVM}}
We now turn our attention to {\tt NSSVM} in \cref{alg:NSSVM}. It reduces to \cref{alg:snsvm} if we set $\sigma=1$, namely, $s_k=s_0$ remains invariant. However, to distinguish them, in the following, we always set $\sigma=1.1$. In this subsection, we focus on how the sparsity level tuning strategy in {\tt NSSVM} affects the final results.  For \cref{alg:snsvm}, we implement it with different fixed $si=s(\beta_i),i=1,2,3$ and name it  {\tt Alg1-si}. In particular, we use $\beta_1=0.4,\beta_2=0.5,\beta_3=0.8$. Recall that {\tt NSSVM}  also requires an initial sparsity level $s_0$ that is set as $s_0=s(0.4)=s1$. Therefore, {\tt NSSVM} and {\tt Alg1-s1} start  with  the  same initial sparsity level.

 \begin{figure}[!th]  
\centering
\begin{subfigure}{.24\textwidth}
	\centering
	\includegraphics[width=1.04\linewidth]{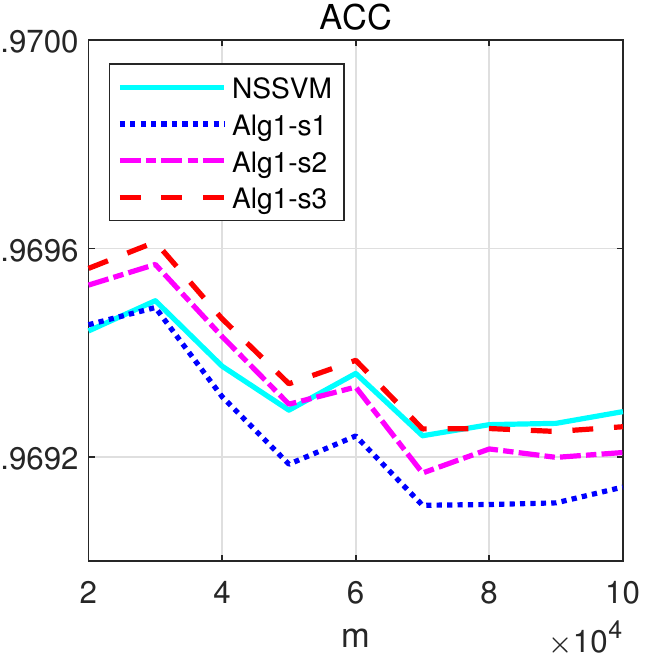}
	%\caption{$m=200$}
	%\label{fig:USairport2010_result2}
\end{subfigure}
\begin{subfigure}{.24\textwidth}
	\centering
	\includegraphics[width=1.04\linewidth]{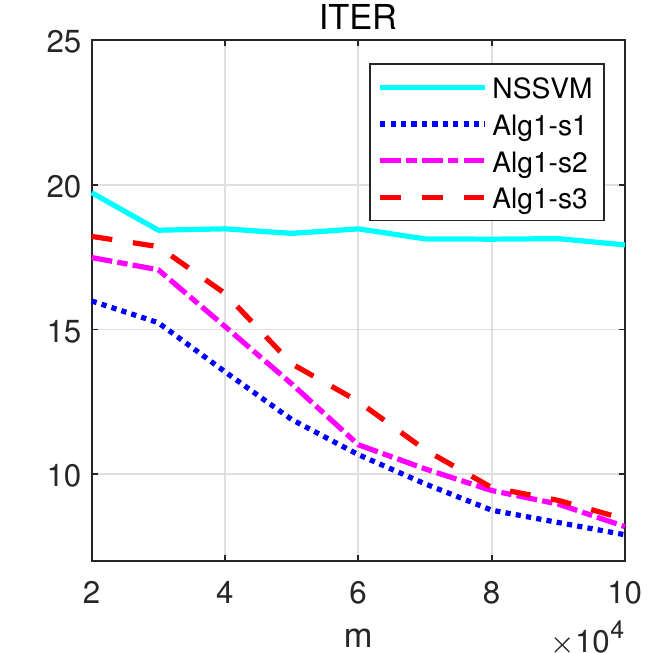}
	%\caption{$m=400$}
	%  \label{fig:Blogs_result2}
\end{subfigure} \\\vspace{1mm}
\begin{subfigure}{.24\textwidth}
	\centering
	\includegraphics[width=1.04\linewidth]{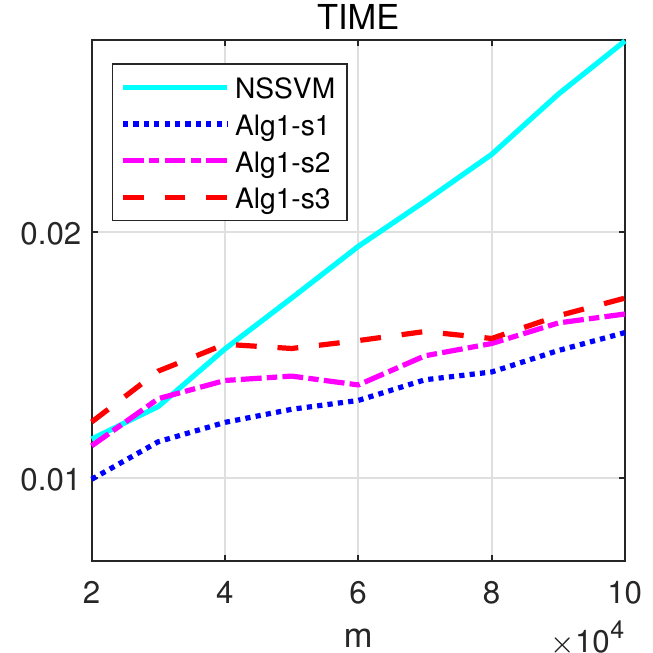}
	%\caption{$m=200$}
	%\label{fig:USairport2010_result2}
\end{subfigure}
\begin{subfigure}{.242\textwidth}
	\centering
	\includegraphics[width=1.04\linewidth]{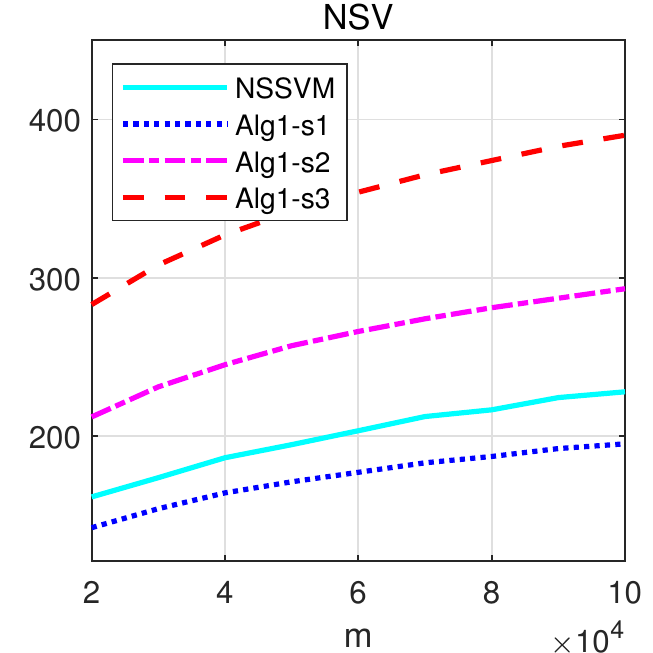}
	%\caption{$m=400$}
	%  \label{fig:Blogs_result2}
\end{subfigure} \\\vspace{-1mm}
\caption{\cref{alg:snsvm} v.s. \cref{alg:NSSVM}}
\label{fig:com-a1-a2}
\end{figure}
The average results of two algorithms  used to solve \cref{ex:syn-data}
 were displayed in Fig. \ref{fig:com-a1-a2}. First, despite  starting  with  the  same initial sparsity level,  {\tt NSSVM}  generated higher {\tt ACC} than {\tt Alg1-s1} because it used more support vectors (i.e., higher {\tt NSV}) due to the tuning strategy.  Second,  {\tt NSSVM} obtained higher  {\tt ACC} than that produced by {\tt Alg1-s2} and  {\tt Alg1-s3} when $m\geq 70000$, but yielded much lower {\tt NSV}.
 This indicated that the tuning strategy performed well to find a  proper sparsity level to maximize the accuracy.  Clearly,  {\tt NSSVM} required longer  time because the tuning strategy led to more steps to meet the stopping criteria \eqref{halt-cods}. Nevertheless, {\tt NSSVM} still ran quickly, requiring less than 0.03 second and  20 steps to converge.

We also compared \cref{alg:snsvm} and {\tt NSSVM} for solving the real datasets, and omitted similar observations. In a nutshell, the tuning strategy enables {\tt NSSVM} to deliver higher accuracy with slightly longer computational time. Based on this,  we next conduct some numerical comparisons between {\tt NSSVM} and some state-of-the-art solvers.

\subsection{Numerical Comparisons}
We note that {\tt NSSVM} involves four parameters $(\eta,C,c,s0)$. On the one hand, it is impractical to derive a universal combination of these four parameters for all datasets. On the other hand, tuning them for each dataset will increase the computational cost. Therefore, in the following numerical experiments, we seek to avoid using different parameters for different datasets. Typically, we fix $\eta=1/m$ and $c=0.01C$, but select  $C$ and $s_0$    as specified in \cref{Table-par}.%$C=2^{-2}$ and $s_0=s(0.05)$ if $m\leq5.8e5$ and $C=\log_2(m)$ and $s_0=s(10)$ otherwise. %We set bigger $\beta$ in $s_0=s(\beta)$   for  {\tt higg}, {\tt susy}, {\tt covt} and {\tt hemp} since their sample sizes are relatively large, comparing with the other datasets. 

\begin{table}[!th] 
	\renewcommand{\arraystretch}{1}\addtolength{\tabcolsep}{-1pt}
	\caption{Selection of parameters.}\vspace{-6mm}
	\label{Table-par}
	\begin{center}
		\begin{tabular}{lllll}  \\	\hline
&cases &  	$C$ &   $s_0$	\\	\hline
\multirow{2}{*}{\cref{ex:syn-data}}&$m\leq10^4$&0.25&$s(0.5)$\\ 
 &$m>10^4$&0.25&$s(1)$\\\hline
&{\tt higg}, {\tt susy}, {\tt covt}, {\tt hemp}	& $\log_{2}(m)$& $s(10)$	\\
\cref{ex:real-data}&{\tt a5a}-{\tt a9a}	& $0.25$& $s(0.2)$	\\
&other	& $0.25$& $s(0.05)$	\\
 \hline
		\end{tabular}
	\end{center}
\end{table}

\subsubsection{Benchmark methods} There is an impressive body of work that has developed numerical methods to tackle the SVM, such as those in \cite{ CL01, PSG02, LL03,  WL07, YL19}. They perform very well, especially for datasets in small or mediate size. Unfortunately, it is unlikely to find a MATLAB implementation for sparse SVM on which we could successfully perform the experiments. For instance, the methods in \cite{LL03} for RSVM \cite{LM01} and the one in \cite{CSS13} were programmed by C++ and couldn't run by MATLAB.
 \begin{table}[!th] 
	\renewcommand{\arraystretch}{1}\addtolength{\tabcolsep}{-1pt}
	\caption{Benchmark methods.}\vspace{-6mm}
	\label{Table-algs}
	\begin{center}
		\begin{tabular}{lllll}  \\	\hline
Algs.&  Source of Code & Ref.& $\ell$   	\\	\hline 
{\tt HSVM} & {\tt libsvm}\footnote{\url{https://www.csie.ntu.edu.tw/~cjlin/libsvm/}} & \cite{CL01} & Hinge loss	\\
 {\tt SSVM}  & {\tt liblssvm}\footnote{\url{https://www.esat.kuleuven.be/sista/lssvmlab/}} & \cite{PSG02}  & Squared Hinge loss	\\ 
 {\tt LSVM} & {\tt liblinear} \footnote{\url{https://www.csie.ntu.edu.tw/~cjlin/liblinear/}} &\cite{FCH08} & Hinge loss\\
   
   {\tt FSVM}  & {\tt fitclinear} \footnote{\url{https://mathworks.com/help/stats/fitclinear.html}} &&  Hinge loss	\\ 
   
    {\tt SNSVM}  &  {\tt Li's lab}\footnote{\url{https://www.researchgate.net/publication/343963444}} & \cite{YL19}& Squared Hinge loss   
\\ \hline
		\end{tabular}
	\end{center}
\end{table}
 
Therefore, we select five solvers in \cref{Table-algs} that do not aim at the sparse SVM. They address the soft-margin SVM models \eqref{SM-SVM} with different loss functions $\ell$. Note that {\tt fitclinear} is a Matlab built-in solver with initializing {\tt Solver} = {\tt `dual'} in order to obtain the number of the support vectors. For the same reason, we set {\tt -s 3} for {\tt FSVM}.  It runs much faster if we set  {\tt -s 2}, but such a setting does not offer the number of support vectors. The other three methods also involve some parameters. For simplicity, all of these parameters are set to their default values.

\subsubsection{Comparisons for datasets with small sizes} 

\begin{figure}[H] 
\centering
\begin{subfigure}{.24\textwidth}
	\centering
	\includegraphics[width=1\linewidth]{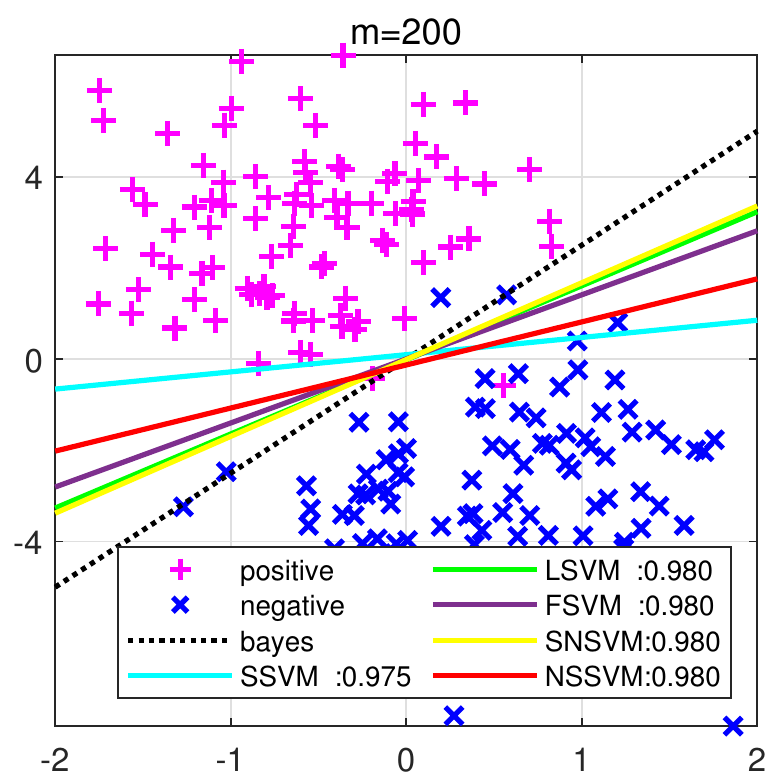}
	%\caption{$m=200$}
	%\label{fig:USairport2010_result2}
\end{subfigure}
\begin{subfigure}{.24\textwidth}
	\centering
	\includegraphics[width=1\linewidth]{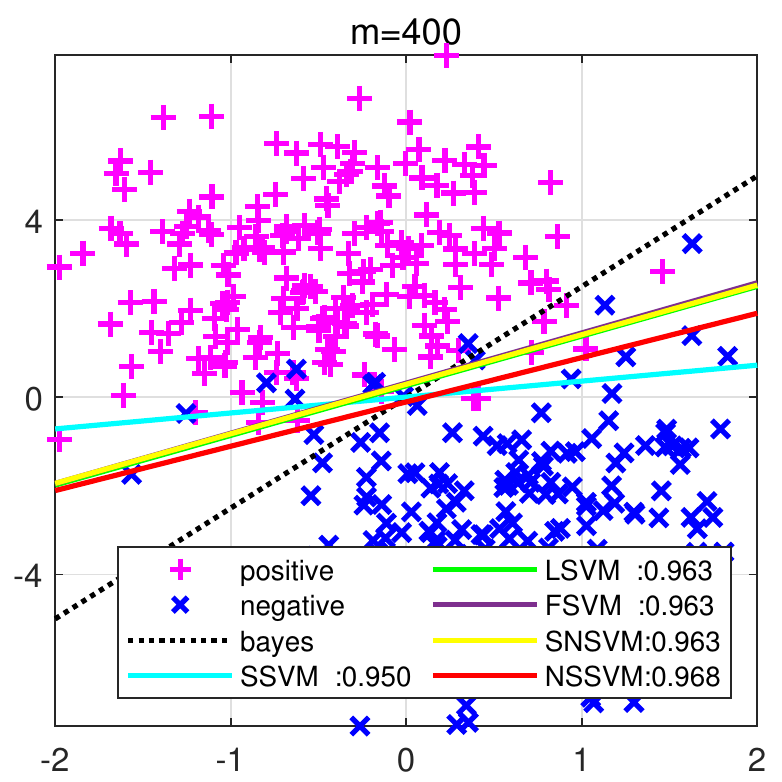}
	%\caption{$m=400$}
	%  \label{fig:Blogs_result2}
\end{subfigure} 
%\\\vspace{1mm}
%\begin{subfigure}{.24\textwidth}
%	\centering
%	\includegraphics[width=1\linewidth]{m-300-eps-converted-to.pdf}
%	%\caption{$m=200$}
%	%\label{fig:USairport2010_result2}
%\end{subfigure}
%\begin{subfigure}{.24\textwidth}
%	\centering
%	\includegraphics[width=1\linewidth]{m-400-eps-converted-to.pdf}
%	%\caption{$m=400$}
%	%  \label{fig:Blogs_result2}
%\end{subfigure} \\
\caption{Classifiers by five solvers for  \cref{ex:syn-data}.}
\label{fig:ex1}
\end{figure}

We first applied five solvers for solving  \cref{ex:syn-data} with small scales $m=200,400$ and depicted the classifiers in  Fig. \ref{fig:ex1} where the Bayes classifier was $w_2=0.25w_1$, see the black dotted line. We did not run
{\tt HSVM}  because it solved the dual problem and did not provide the solution $(\bfw, b)$.  Overall, all solvers can classify the dataset well.

\begin{figure}[!th]  
\centering
\begin{subfigure}{.24\textwidth}
	\centering
	\includegraphics[width=1.05\linewidth]{ 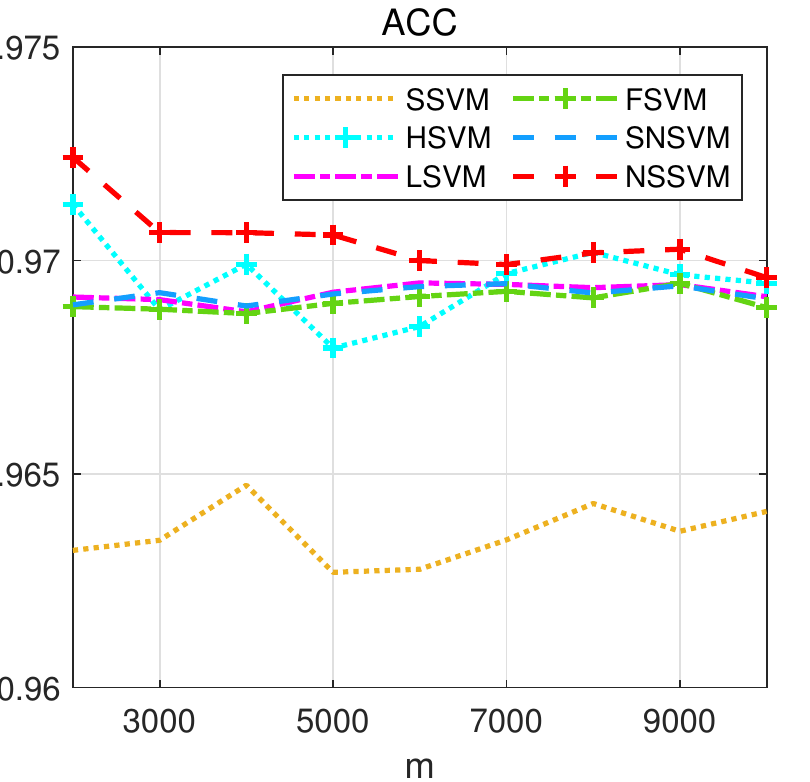}
	%\caption{$m=200$}
	%\label{fig:USairport2010_result2}
\end{subfigure}
\begin{subfigure}{.24\textwidth}
	\centering
	\includegraphics[width=1.05\linewidth]{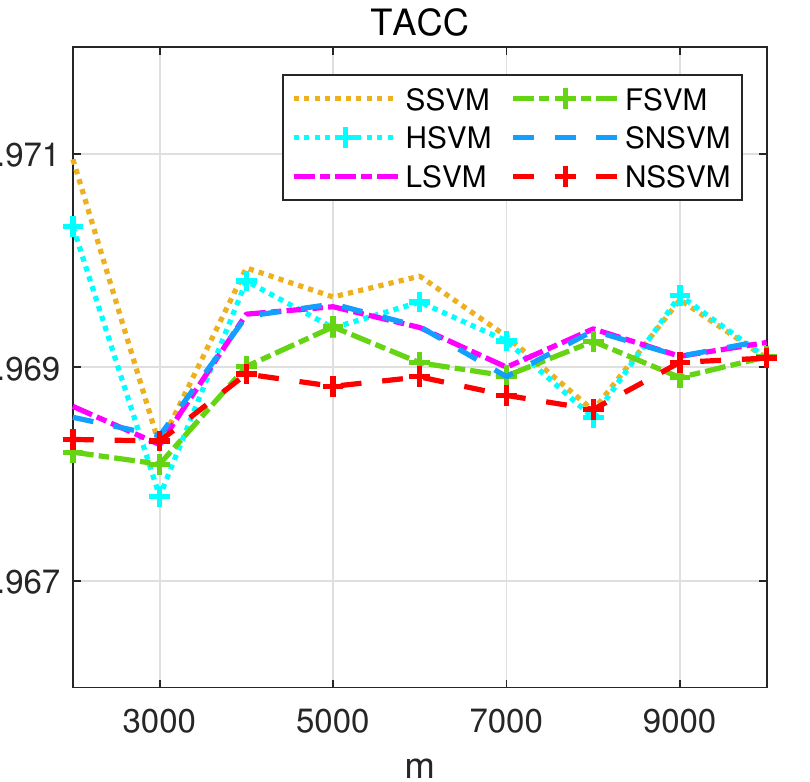}
	%\caption{$m=400$}
	%  \label{fig:Blogs_result2}
\end{subfigure} \\\vspace{1mm}
\begin{subfigure}{.24\textwidth}
	\centering
	\includegraphics[width=1.05\linewidth]{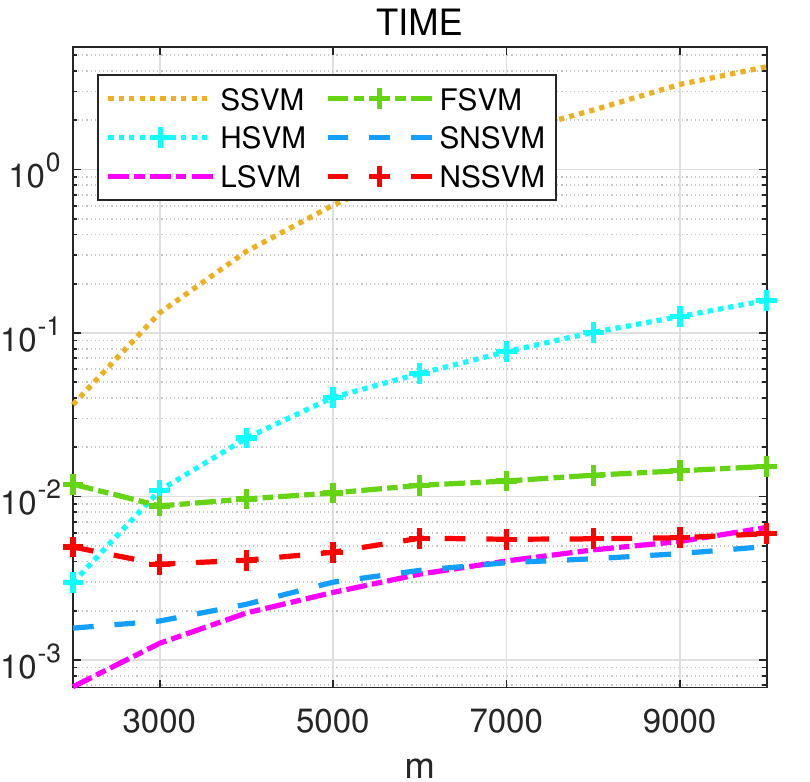}
	%\caption{$m=200$}
	%\label{fig:USairport2010_result2}
\end{subfigure}
\begin{subfigure}{.242\textwidth}
	\centering
	\includegraphics[width=1.05\linewidth]{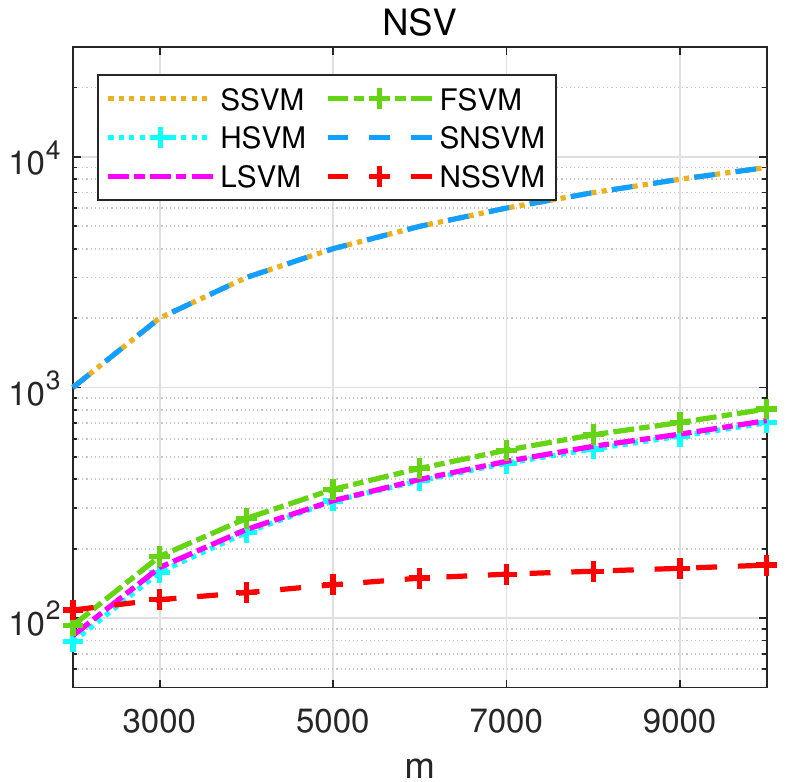}
	%\caption{$m=400$}
	%  \label{fig:Blogs_result2}
\end{subfigure} \\
\caption{Average results of six solvers for   \cref{ex:syn-data}.}
\label{fig:eff-6alg-m} 
\end{figure}

  We then used six solvers to address \cref{ex:syn-data} with different sample sizes $m\in[2000,10^4]$.  Average results over 100 trials were reported in Fig. \ref{fig:eff-6alg-m}. It is observed that {\tt NSSVM} obtained the largest {\tt ACC}, but rendered smallest {\tt TACC} for most cases. This was probably due to  the smallest number of support vectors used.  Evidently, {\tt LSVM} and {\tt SNSVM} ran the fastest, followed by {\tt FSVM} and {\tt NSSVM}.

 Next,  six methods were used to handle   \cref{ex:real-data} with small scale datasets. As shown in \cref{table:ex2-small}, where `$--$' denoted that  {\tt SSVM} required memory that was beyond the capacity of our laptop; here,  again {\tt NSSVM} did not show the advantage of delivering the accuracy. This was not surprising because it used a relatively small number of support vectors. However, it ran the fastest for most datasets.
 
To summarize,  {\tt NSSVM} did not present the advantage of higher accuracy for dealing with small scale datasets. One reason for this is that the solution to the model \eqref{SM-SVM-h-equ-dual-sparse} may not be sufficiently sparse for small scale datasets. In other words, there is no strong need to perform data reduction for these datasets. Hence, in the sequel, we employ four methods to address test problems with datasets on much larger scales. However, we remove {\tt HSVM} and {\tt SSVM} from comparisons in the following since they either are too slow or required memory beyond the capability of our device.

\begin{table}[!th]
\centering
\caption{Results of six solvers  solving   \cref{ex:real-data} with datasets with small sizes.}
\label{table:ex2-small} 
\renewcommand{\arraystretch}{1}\addtolength{\tabcolsep}{1pt}	
\begin{tabular}{lcccccc}\hline

		&	{\tt SSVM}	&	{\tt HSVM}	&	{\tt LSVM}	&	{\tt FSVM}	&	{\tt SNSVM}	&	{\tt NSSVM}		\\ \cline{2-7} 
		&\multicolumn{6}{c}{{\tt ACC} (\%)} \\
%  \multicolumn{3}{c}{{\tt TIME} (in seconds)}&& \multicolumn{3}{c}{{\tt NSV}} \\ 
 \cline{2-7} 
{\tt	pcb1}&	96.21 	&	96.22 	&	96.22 	&	96.22 	&	96.24 	&	96.24 	\\
{\tt	pcb2}&	95.95 	&	95.96 	&	95.96 	&	95.96 	&	{\bf 95.97} 	&	95.95 	\\
{\tt	pcb3}&	95.31 	&	95.31 	&	95.31 	&	95.31 	&	95.30 	&	95.31 	\\
{\tt	pcb4}&	94.71 	&	94.74 	&	94.74 	&	94.74 	&	94.72 	&	94.71 	\\
{\tt	pcb5}&	93.04 	&	93.12 	&	93.16 	&	93.16 	&	{\bf 93.18} 	&	93.04 	\\
{\tt	a5a}&	24.46 	&	85.06 	&	85.05 	&	81.82 	&	{\bf 85.16} 	&	84.10 	\\
{\tt	a6a}&	$--$	&	{\bf 84.89} 	&	84.87 	&	76.11 	&	84.80 	&	83.95 	\\
{\tt	a7a}&	$--$	&	84.88 	&	{\bf 84.96} 	&	81.93 	&	84.89 	&	84.16 	\\
{\tt	a8a}&	$--$	&	84.69 	&	84.67 	&	83.00 	&	{\bf 84.74} 	&	84.01 	\\
{\tt	a9a}&	$--$	&	84.99 	&	84.99 	&	80.63 	&	84.96 	&	84.04 	\\\hline
&\multicolumn{6}{c}{{\tt TACC} (\%)} \\
 \cline{2-7} 	\hline
{\tt	pcb}&	{\bf 95.73} 	&	95.58 	&	95.58 	&	95.58 	&	95.58 	&	95.58 	\\
{\tt	pcb2}&	96.85 	&	97.15 	&	97.15 	&	97.15 	&	97.15 	&	97.15 	\\
{\tt	pcb3}&	95.05 	&	95.05 	&	95.05 	&	95.05 	&	95.05 	&	95.05 	\\
{\tt	pcb4}&	94.69 	&	94.99 	&	94.99 	&	94.99 	&	94.89 	&	94.99 	\\
{\tt	pcb5}&	93.06 	&	93.23 	&	93.23 	&	93.23 	&	93.23 	&	93.23 	\\
{\tt	a5a}&	 84.42 	&	84.39 	&	84.39 	&	81.07 	&	{\bf 84.50} 	&	83.60 	\\
{\tt	a6a}&	$--$	&	84.72 	&	84.66 	&	76.12 	&	{\bf 84.77} 	&	84.13 	\\
{\tt	a7a}&	$--$	&	84.84 	&	84.82 	&	82.26 	&	{\bf 85.00} 	&	84.00 	\\
{\tt	a8a}&	$--$	&	85.16 	&	85.13 	&	83.47 	&	{\bf 85.42} 	&	84.59 	\\
{\tt	a9a}&	$--$	&	  84.98 	&	{\bf 85.00} 	&	80.96 	&	84.94 	&	84.47 	\\\hline
	
&  \multicolumn{6}{c}{{\tt TIME} (in seconds)} \\ 
 \cline{2-7}
{\tt pcb1}&	1.905 	&	0.702 	&	{\bf 0.037} 	&	0.241 	&	0.057 	&	0.082 	\\
{\tt	pcb2}&	4.656 	&	1.577 	&	0.051 	&	0.243 	&	0.064 	&	{\bf 0.039} 	\\
{\tt	pcb3}&	4.987 	&	1.761 	&	0.046 	&	0.064 	&	0.044 	&	{\bf 0.019} 	\\
{\tt	pcb4}&	4.282 	&	1.722 	&	0.048 	&	0.051 	&	0.038 	&	{\bf 0.016} 	\\
{\tt	pcb5}&	1.241 	&	0.805 	&	0.044 	&	0.033 	&	0.021 	&	{\bf 0.013} 	\\
{\tt	a5	}&	29.72 	&	5.414 	&	{\bf 0.030} 	&	0.110 	&	0.044 	&	0.049 	\\
{\tt	a6a}&	$--$	&	11.15 	&	{\bf 0.051} 	&	0.247 	&	0.107 	&	0.085  	\\
{\tt	a7a}&	$--$	&	19.15 	&	0.067 	&	{\bf 0.072} 	&	0.179 	&	0.097 	\\
{\tt	a8a}&	$--$	&	33.07 	&	0.105 	&	{\bf 0.067} 	&	0.287 	&	0.121 	\\
{\tt	a9a}&	$--$	&	70.59 	&	0.158 	&	{\bf 0.072} 	&	0.361 	&	0.130 	\\
 
 \hline
 &  \multicolumn{6}{c}{{\tt NSV}} \\ 
 \cline{2-7}
{\tt	pcb1}&	6325	&	497	&	2320	&	2263	&	6325	&	{\bf 190}	\\
{\tt	pcb2}&	9156	&	791	&	4166	&	3557	&	9156	&	{\bf 182}	\\
{\tt	pcb3}&	9453	&	981	&	4692	&	4375	&	9453	&	{\bf 184}	\\
{\tt	pcb4}&	8813	&	990	&	4538	&	4762	&	8813	&	{\bf 179}	\\
{\tt	pcb5}&	5319	&	758	&	2167	&	2440	&	5319	&	{\bf 145}	\\
{\tt	a5a}&	6414	&	2294	&	2332	&	3891	&	6414	&	{\bf 882}	\\
{\tt	a6a}&	$--$	&	4011	&	4082	&	6910	&	11220	&	{\bf 1148}	\\
{\tt	a7a}&	$--$	&	5757	&	5899	&	9931	&	16100	&	{\bf 1339	}\\
{\tt	a8a}&	$--$	&	8135	&	8292	&	13992	&	22696	&	{\bf 1534}	\\
{\tt	a9a}&	$--$	&	11533	&	11772	&	20030	&	32561	&	{\bf 1754}	\\
 
 \hline
\end{tabular}
\end{table}

\subsubsection{Comparisons for datasets with large sizes} 
This part focus on instances with much larger scales. First, we employ four algorithms to solve \cref{ex:syn-data} with  $m$ from $[10^5, 10^8]$.  The data reported in \cref{table:ex1-less-n-bigger-m-2} were the results averaged  over 20 trials.  Since {\tt LSVM} ran for too long when $m>10^6$, its results were omitted. The training  classification accuracy values for {\tt NSSVM}, {\tt SNSVM} and {\tt LSVM}  were similar and were better than those  obtained by {\tt FSVM}. However, it is clearly observed that {\tt NSSVM} ran the fastest and used much fewer support vectors.  Moreover, the superiority of {\tt NSSVM} becomes more evident with larger $m$.

\begin{table}[H] 
\centering
\caption{Results of four solvers  solving   \cref{ex:syn-data} with large-scale datasets.}
\label{table:ex1-less-n-bigger-m-2} 
\renewcommand{\arraystretch}{1}\addtolength{\tabcolsep}{-3.2pt}	
\begin{tabular}{lcccc c cccc}\hline
	$m$&  	{\tt LSVM} &	  {\tt FSVM} &	{\tt SNSVM} &	{\tt NSSVM}	&&
	      	{\tt LSVM} &	  {\tt FSVM} &	{\tt SNSVM} &	{\tt NSSVM}\\ \hline
	      	&\multicolumn{4}{c}{{\tt ACC} (\%)}&& \multicolumn{4}{c}{{\tt TACC} (\%)} \\ 
 \cline{2-5}\cline{7-10}
%$10^4$	&	96.95 	&	96.94 	&	96.95 	&	96.95 	&	&	96.97 	&	96.95 	&	96.97 	&	96.96 	\\
$10^5$	&	96.95 	&	96.91 	&	96.95 	&	96.95 	&	&	96.96 	&	96.93 	&	96.96 	&	96.94 	\\
$10^6$	&	96.94 	&	96.91 	&	96.94 	&	96.94 	&	&	96.94 	&	96.93 	&	96.94 	&	96.93 	\\
$10^7$	&	$--$	&	96.92 	&	96.93 	&	96.93 	&	&	$--$	&	96.93 	&	96.93 	&	96.93 	\\
$10^8$	&	$--$	&	96.91 	&	96.93 	&	96.93 	&	&	$--$	&	96.93 	&	96.93 	&	96.93 	\\

 \hline
&
  \multicolumn{4}{c}{{\tt TIME} (in seconds)}&& \multicolumn{4}{c}{{\tt NSV}$/m$}
\\ \cline{2-5}\cline{7-10}
%$10^4$	&	0.01 	&	0.04 	&	0.01 	&	0.01 	&	&	7.78e-2	&	8.79e-2	&	1	&	3.33e-2	\\
$10^5$	&	0.10 	&	0.10 	&	0.04 	&	0.04 	&	&	7.83e-2	&	8.85e-2	&	1	&	{\bf 5.94e-3}	\\
$10^6$	&	3.88 	&	2.43 	&	0.36 	&	{\bf 0.33} 	&	&	7.81e-2	&	8.87e-2	&	1	&	{\bf 8.62e-4}	\\
$10^7$	&	$--$	&	30.62 	&	3.42 	&	{\bf 2.08} 	&	&	$--$	&	8.88e-2	&	1	&	{\bf 1.09e-4}	\\
$10^8$	&	$--$	&	447.1 	&	34.1 	&	{\bf 12.3} 	&	&	$--$	&	8.89e-2	&	1	&	{\bf 1.44e-5}	\\
\hline
\end{tabular}
\end{table}

\begin{table}[!th] 
\centering
\caption{Results of four solvers solving \cref{ex:real-data} with datasets with large sizes.}
\label{table:ex2-less-n-bigger-m-3} 
\renewcommand{\arraystretch}{1}\addtolength{\tabcolsep}{-4pt}	
\begin{tabular}{lcccc c cccc }\hline
			&	{\tt LSVM}	&	{\tt FSVM}	&	{\tt SNSVM}	&	{\tt NSSVM}	&&	{\tt LSVM}	&	{\tt FSVM}	&	{\tt SNSVM}	&	{\tt NSSVM}	\\\cline{2-5}\cline{7-10} 
			&\multicolumn{4}{c}{{\tt ACC} (\%)}&& \multicolumn{4}{c}{{\tt TACC} (\%)}\\\cline{2-5}\cline{7-10} 
{\tt	 mrpe 	}	&	95.01 	&	94.78 	&	{\bf 95.08} 	&	95.00 	&&	95.25 	&	94.90 	&	95.25 	&	95.25 	\\
{\tt	 dhrb 	}	&	82.70 	&	82.71 	&	{\bf 83.25} 	&	82.76 	&&	83.43 	&	83.43 	&	{\bf 84.11} 	&	 83.48 	\\
{\tt	 aips 	}	&	{\bf 87.66} 	&	87.26 	&	87.42 	&	86.54 	&&	{\bf 87.41} 	&	87.08 	&	87.13 	&	85.99 	\\
{\tt	 sctp 	}	&	89.95 	&	78.01 	&	{\bf 91.20} 	&	90.50 	&&	90.00 	&	77.80 	&	{\bf 91.19} 	&	90.51 	\\
{\tt	 skin 	}	&	{\bf 92.90} 	&	 92.73 	&	92.37 	&	91.04 	&&	{\bf 92.66} 	&	92.45 	&	92.23 	&	90.63 	\\
{\tt	 ccfd 	}	&	99.94 	&	99.94 	&	99.92 	&	99.94 	&&	99.92 	&	99.92 	&	99.88 	&	99.92 	\\
{\tt	 rlc1 	}	&	100.0 	&	100.0 	&	100.0 	&	100.0 	&&	100.0 	&	100.0 	&	100.0 	&	100.0 	\\
{\tt	 rlc2 	}	&	100.0 	&	100.0 	&	100.0 	&	100.0 	&&	100.0 	&	100.0 	&	100.0 	&	100.0 	\\
{\tt	 rlc3 	}	&	100.0 	&	100.0 	&	100.0 	&	100.0 	&&	100.0 	&	100.0 	&	99.99 	&	100.0 	\\
{\tt	 rlc4 	}	&	100.0 	&	100.0 	&	100.0 	&	100.0 	&&	100.0 	&	100.0 	&	100.0 	&	100.0 	\\
{\tt	 rlc5 	}	&	100.0 	&	100.0 	&	100.0 	&	100.0 	&&	100.0 	&	100.0 	&	100.0 	&	100.0 	\\
{\tt	 rlc6 	}	&	100.0 	&	100.0 	&	100.0 	&	100.0 	&&	100.0 	&	100.0 	&	99.99 	&	99.99 	\\
{\tt	 rlc7 	}	&	100.0 	&	100.0 	&	100.0 	&	100.0 	&&	100.0 	&	100.0 	&	100.0 	&	100.0 	\\
{\tt	 rlc8 	}	&	100.0 	&	100.0 	&	100.0 	&	100.0 	&&	100.0 	&	100.0 	&	100.0 	&	100.0 	\\
{\tt	 rlc9 	}	&	100.0 	&	100.0 	&	100.0 	&	100.0 	&&	100.0 	&	100.0 	&	100.0 	&	100.0 	\\
{\tt	 rlc10 	}	&	100.0 	&	100.0 	&	100.0 	&	100.0 	&&	100.0 	&	100.0 	&	100.0 	&	100.0 	\\
{\tt	 covt 	}	&	{\bf 76.29} 	&	75.73 	&	75.67 	&	75.67 	&&	{\bf 76.15} 	&	75.47 	&	75.54 	&	75.51 	\\
{\tt	 susy 	}	&	78.44 	&	78.55 	&	78.66 	&	{\bf 78.82} 	&&	78.39 	&	78.52 	&	78.63 	&	{\bf 78.79} 	\\
{\tt	 hepm 	}	&	78.36 	&	83.31 	&	83.63 	&	{\bf 83.74} 	&&	78.33 	&	83.23 	&	83.61 	&	{\bf 83.71} 	\\
{\tt	 higg 	}	&	47.01 	&	63.84 	&	64.13 	&	{\bf 64.16} 	&&	46.98 	&	63.79 	&	64.07 	&	{\bf 64.08} 	\\\hline

  &\multicolumn{4}{c}{{\tt TIME} (in seconds)}&& \multicolumn{4}{c}{{\tt NSV}$/m$} \\ 	\cline{2-5}\cline{7-10} 	
{\tt	 mrpe 	}	&	31.75 	&	2.28 	&	5.37 	&	{\bf 1.06} 	&&	8.23e-1	&	5.60e-1	&	1	&	{\bf 3.63e-2}	\\
{\tt	 dhrb 	}	&	0.17 	&	0.10 	&	0.07 	&	{\bf 0.04} 	&&	7.68e-1	&	6.92e-1	&	1	&	{\bf 2.59e-3}	\\
{\tt	 aips 	}	&	0.87 	&	0.21 	&	0.18 	&	{\bf 0.09} 	&&	3.01e-1	&	5.30e-1	&	1	&	{\bf 2.33e-3}	\\
{\tt	 sctp 	}	&	13.48 	&	1.09 	&	10.61 	&	{\bf 0.38} 	&&	3.82e-1	&	5.79e-1	&	1	&	{\bf 5.84e-3}	\\
{\tt	 skin 	}	&	0.42 	&	0.33 	&	0.11 	&	{\bf 0.08} 	&&	1.97e-1	&	2.55e-1	&	1	&	{\bf 2.20e-4}	\\
{\tt	 ccfd 	}	&	2.66 	&	0.52 	&	0.40 	&	{\bf 0.16} 	&&	5.56e-3	&	1.01e-2	&	1	&	{\bf 9.41e-4}	\\
{\tt	 rlc1 	}	&	0.39 	&	0.66 	&	0.25 	&	{\bf 0.11} 	&&	3.32e-4	&	5.95e-4	&	1	&	{\bf 2.17e-4}	\\
{\tt	 rlc2 	}	&	0.39 	&	0.51 	&	0.26 	&	{\bf 0.17} 	&&	2.47e-4	&	5.13e-4	&	1	&	{\bf 2.40e-4}	\\
{\tt	 rlc3 	}	&	0.49 	&	0.56 	&	0.23 	&	{\bf 0.12} 	&&	2.85e-4	&	5.86e-4	&	1	&	{\bf 2.17e-4}	\\
{\tt	 rlc4 	}	&	0.38 	&	0.67 	&	0.25 	&	{\bf 0.18} 	&&	2.68e-4	&	5.79e-4	&	1	&	{\bf 2.40e-4}	\\
{\tt	 rlc5 	}	&	0.40 	&	0.57 	&	0.26 	&	{\bf 0.12} 	&&	2.89e-4	&	5.76e-4	&	1	&	{\bf 2.17e-4}	\\
{\tt	 rlc6 	}	&	0.50 	&	0.51 	&	0.27 	&	{\bf 0.2}1 	&&	2.44e-4	&	5.51e-4	&	1	&	{\bf 2.64e-4}	\\
{\tt	 rlc7 	}	&	0.38 	&	0.52 	&	0.25 	&	{\bf 0.12} 	&&	2.84e-4	&	5.31e-4	&	1	&	{\bf 2.17e-4}	\\
{\tt	 rlc8 	}	&	0.68 	&	0.67 	&	0.25 	&	{\bf 0.12} 	&&	2.68e-4	&	5.77e-4	&	1	&	{\bf 2.17e-4}	\\
{\tt	 rlc9 	}	&	0.66 	&	0.54 	&	0.25 	&	{\bf 0.12} 	&&	3.04e-4	&	6.00e-4	&	1	&	{\bf 2.17e-4}	\\
{\tt	 rlc10 	}	&	0.39 	&	0.55 	&	0.25 	&	{\bf 0.12} 	&&	3.25e-4	&	5.91e-4	&	1	&	{\bf 2.17e-4}	\\
{\tt	 covt 	}	&	13.23 	&	2.05 	&	9.22 	&	{\bf 1.11} 	&&	4.75e-1	&	5.93e-1	&	1	&	{\bf 1.61e-1}	\\
{\tt	 susy 	}	&	596.3	&	16.69 	&	16.44 	&	{\bf 2.43} 	&&	4.15e-1	&	4.74e-1	&	1	&	{\bf 1.26e-2}	\\
{\tt	 hepm 	}	&	1688 	&	28.39 	&	26.97 	&	{\bf 3.62} 	&&	3.45e-1	&	4.76e-1	&	1	&	{\bf 1.27e-2}	\\
{\tt	 higg 	}	&	2938 	&	46.86 	&	61.75 	&	{\bf 4.78} 	&&	2.37e-1	&	7.30e-1	&	1	&	{\bf 8.56e-3}	\\

 \hline
\end{tabular}
\end{table}

Results obtained by four solvers for classifying real datasets in \cref{ex:real-data} on larger scales were reported in   \cref{table:ex2-less-n-bigger-m-3}. For the accuracy, there was no significant difference of  {\tt ACC} and {\tt TACC} produced by four solvers. For the computational speed, it is clearly observed that  {\tt NSSVM} ran the fastest for all datasets. With regard to the number of support vectors,  {\tt NSSVM} outperformed the other methods because it obtained the smallest {\tt NSV}. Taking {\tt higg} as an instance, the other three solvers respectively produced $0.237m, 0.73m$,  and $m$ support vectors and consumed $2938$, $46.86$ and $61.75$ seconds to classify the data. By contrast, our proposed method used $0.00856m$ support vectors and only required $4.78$ seconds to achieve the highest classification accuracy.

In summary,   {\tt NSSVM} delivered the desirable classification accuracy and ran rapidly by using a relatively small number of support vectors for large-scale datasets.

\section{Conclusion}
 {The sparsity constraint was employed in the kernel-based SVM model \eqref{SM-SVM-h-equ-dual-sparse} allowing for the control of the number of support vectors. Therefore, the  memory required for data storage and the computational cost were significantly reduced, enabling large scale computation.  We feel that the established theory  and the proposed Newton type method deserve to be further explored for using the SVM model with some nonlinear kernel-based models.}

\section*{Acknowledgements}
I sincerely thank the two referees and Professor Naihua Xiu from Beijing Jiaotong University  for their valuable suggestions to improve this paper.

\ifCLASSOPTIONcompsoc


\ifCLASSOPTIONcaptionsoff
  \newpage
\fi



%\bibliographystyle{IEEEtran}

\begin{IEEEbiography}[{\includegraphics[width=1in,clip,keepaspectratio]{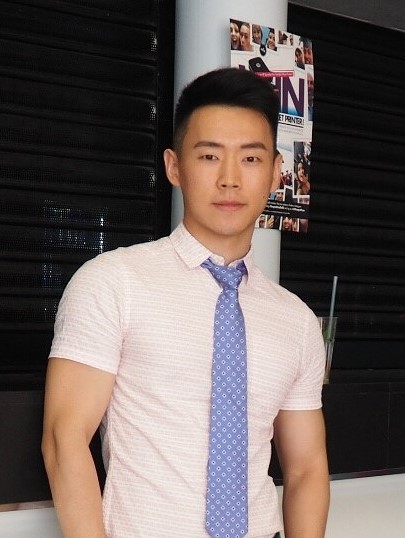}}]{Shenglong Zhou} received the B.Sc. degree in
information and computing science in 2011 and
the M.Sc. degree in operational research in 2014
from Beijing Jiaotong University, China, and the
Ph.D. degree in operational research in 2018
from the University of Southampton, the United
Kingdom, where he was the Research Fellow
from 2017 to 2019 and is currently a Teaching
Fellow. His research interests include the theory
and methods of optimization in the fields of
sparse, low-rank matrix and bilevel optimization.
\end{IEEEbiography}


%\newpage


%\vfill


\appendices
\section{Proofs of all theorems}

\subsection*{Proof  of \cref{the:duality}}
%\textcolor{blue}
{{\bf proof} By introducing $u_i=1-y_{i} (\langle \bfw, \bfx_i\rangle +b)$, the problem \eqref{SM-SVM-cc} is rewritten as
\begin{eqnarray}
\label{SM-SVM-h-equ}
\min_{\bfw,\bfu}&&\frac{1}{2}\|\bfw\|^2+  \sum_{i=1}^{m}\ell_{cC}\left(u_{i}\right), \\
{\rm s.t.}&&u_i+y_{i} (\langle \bfw, \bfx_i\rangle +b) =1,~ i\in[m].\nonumber
\end{eqnarray} 
The Lagrangian function of the above problem is
\begin{eqnarray*}
L(\bfw,\bfu,b,\bfal)&:=& \frac{1}{2}\|\bfw\|^2+  \sum_{i=1}^{m}\ell_{cC}\left(u_{i}\right)\\
& -&\sum_{i=1}^{m} \alpha_i \Big(u_i+y_{i}( \langle \bfw, \bfx_i\rangle+b)-1\Big). 
\end{eqnarray*} 
The dual problem of \eqref{SM-SVM-h-equ} is then given by
\begin{eqnarray}\label{dual-problem}
\max_\bfal \min_{\bfw,\bfu,b} L(\bfw,\bfu,b,\bfal). 
\end{eqnarray} 
For a fixed $\bfal$, the inner problem $\min_{\bfw,\bfu,b} L(\bfw,\bfu,b,\bfal)$ has three variables that are independent of each other. Hence, it can be divided into three sub-problems
 \allowdisplaybreaks \begin{eqnarray}\label{three-sub}
&&\min_{\bfw,\bfu,b}~~ L(\bfw,\bfu,b,\bfal) \\
&=&\min_{\bfw}~\frac{1}{2}\|\bfw\|^2-\sum_{i=1}^{m} \alpha_i y_{i}  \langle \bfw, \bfx_i\rangle +\sum_{i=1}^{m} \alpha_i  \nonumber\\
&+& \min_{\bfu} 
\sum_{i=1}^{m}\Big(  \ell_{cC}(u_{i})-\alpha_i  u_i\Big)- \underset{b}{\min} \sum_{i=1}^{m}  \alpha_i  y_i b. \nonumber 
\end{eqnarray}
For the sub-problem with respect to $\bfw$, the optimal solution $\bfw$ is attained at $w_j-\sum_{i=1}^{m} \alpha_i y_{i} x_{ij}=0, i\in[n]$, namely, 
\begin{eqnarray}
\label{SM-SVM-h-equ-dual-sub1-opt}
 \bfw=Q\bfal.
\end{eqnarray} 
using notation $Q$ in \cref{notation}. Substituting \eqref{SM-SVM-h-equ-dual-sub1-opt} into the following objective function yields 
\begin{eqnarray}
\label{SM-SVM-h-equ-dual-sub1}
\underset{\bfw}{\min}~\frac{1}{2}\|\bfw\|^2-\sum_{i=1}^{m} \alpha_i y_{i}  \langle \bfw, \bfx_i\rangle=-\frac{1}{2}\| Q \bfal\|^2.
%\nonumber\\
% &=& \frac{1}{2}\|Q\bfal\|^2-\sum_{i=1}^{m} \alpha_i y_{i}  \langle Q\bfal, \bfx_i\rangle\nonumber\\
%&=&\frac{1}{2}\|Q\bfal\|^2-\sum_{i=1}^{m}   \langle Q\bfal, Q\bfal\rangle\nonumber\\
%&=&-\frac{1}{2}\| Q \bfal\|^2.
\end{eqnarray} 
For the sub-problem with respect to $\bfu$, the optimal solution $\bfu$ is attained at, for $i\in[m]$, \begin{eqnarray}
\label{SM-SVM-h-equ-dual-sub2-opt}
0= \ell'_{c,C}(u_{i})- \alpha_i~~\Longleftrightarrow~~\alpha_i=\left\{
\begin{array}{rr}
Cu_i, & u_i\geq0,\\
cu_i,& u_i< 0,
\end{array}
\right.
\end{eqnarray} 
 where $\ell'_{c,C}(u_{i})$ is the derivative of $\ell'_{c,C}$ at $u_i$ that leads to \begin{eqnarray*}
 \ell_{cC}(u_{i})- \alpha_i  u_i&=& 
\begin{cases}
{C}u_i^2/2- \alpha_i  u_i, &u_i\geq0,\\
{c}u_i^2/2- \alpha_i  u_i,& u_i< 0,
\end{cases}\\
&=& 
\begin{cases}
-{\alpha_i^2}/{(2C)}, &\alpha_i\geq0,\\
-{\alpha_i^2}/{(2c)}, &\alpha_i< 0,\\
\end{cases} \\
&=& -h_{cC}(\alpha_i ).
\end{eqnarray*}
Therefore, we have the sub-problem 
\begin{eqnarray}
\label{SM-SVM-h-equ-dual-sub2}
\underset{\bfu}{\min}~
 \sum_{i=1}^{m} \Big( \ell_{cC}(u_{i})- \alpha_i  u_i \Big) = -\sum_{i=1}^{m}h_{cC}(\alpha_i ).
\end{eqnarray} 
For the third sub-problem with respect to $b$, the optimal solution is attained at 
\begin{eqnarray}
\label{SM-SVM-h-equ-dual-sub3-opt}
\langle\bfal,\bfy\rangle= \sum_{i=1}^{m}  \alpha_i  y_i= 0.
\end{eqnarray} 
This results in
 $$ {\min}_b 
 \sum_{i=1}^{m}  \alpha_i  y_i b =0.$$
 Overall, the sub-problem \eqref{three-sub} becomes
 \begin{eqnarray}
\min_{\bfw,\bfu,b}~L(\bfw,\bfu,b,\bfal) = -\frac{1}{2}\| Q \bfal\|^2 +\sum_{i=1}^{m} \alpha_i   -h_{cC}(\alpha_i )
\end{eqnarray}  
 where $\bfal$ satisfies \eqref{SM-SVM-h-equ-dual-sub3-opt}. Hence, the dual problem \eqref{dual-problem} is equivalent to 
 \begin{eqnarray}\label{dual-pro}
&&\min_\bfal\Big\{-\min_{\bfw,\bfu,b} L(\bfw,\bfu,b,\bfal)\Big\}\\
&=&\min_\bfal\left\{\frac{1}{2}\| Q \bfal\|^2 -\sum_{i=1}^{m} \alpha_i + h_{cC}(\alpha_i ):  \langle\bfal,\bfy\rangle=0\right\}, \nonumber
\end{eqnarray} 
which is same as \eqref{SM-SVM-h-equ-dual}. Let $(\widehat\bfw,\widehat\bfu,\widehat b)$ and $\bfal^*$ be the optimal solutions to the primal problem \eqref{SM-SVM-h-equ} and dual problem \eqref{dual-pro}, respectively. Then from \eqref{SM-SVM-h-equ-dual-sub1-opt}, $\widehat\bfw  =Q\bfal^*$.  For $\widehat b$, it follows 
\begin{eqnarray*}
y_{i} \widehat b&\overset{\eqref{SM-SVM-h-equ}}{=}&1-  \widehat u_i-  \langle \widehat\bfw, y_{i}\bfx_i\rangle  \\
&\overset{\eqref{SM-SVM-h-equ-dual-sub1-opt}}{=}&1-  \widehat u_i-  \langle Q\bfal^*, y_{i}\bfx_i\rangle  \\
&\overset{\eqref{SM-SVM-h-equ-dual-sub2-opt},\eqref{diag-hessian}}{=}&1- E_{ii}(\bfal^*)\alpha_i^*-  \langle Q\bfal^*, y_{i}\bfx_i\rangle,\\
&\overset{\eqref{SM-SVM-h-equ-dual-sparse-opt}}{=}&1- E_{iS_*}(\bfal^*)\alpha_{S_*}^*-  \langle Q_{S_*}\bfal^*_{S_*}, Q_{i}\rangle,\\
&\overset{\eqref{gradient-hessian}}{=}&1- H_{iS_*}(\bfal^*)\alpha_{S_*}^*,
\end{eqnarray*} 
for any $i\in S_*=\supp(\bfal^*)$. This leads to
\begin{eqnarray*}
%\widehat b=y_{i}- y_{i}E_{ii}(\bfal^*)\alpha_i^*-  \langle Q\bfal^*, \bfx_i\rangle
 \widehat b= y_i- y_i H_{i S_*}(\bfal^*)\alpha_{S_*}^*, \forall~i\in S_*.
\end{eqnarray*} 
Summing both sides of the above equation for all $i\in S_*$ yields $ \widehat b$ in \eqref{w-b}, completing the proof.  
\hfill   $ \Box $}

\subsection*{Proof of  \cref{existence-global-minimizers}}
{\bf proof} The solution set is nonempty since $0$ satisfies the constraints of \eqref{SM-SVM-h-equ-dual-sparse}.  The problem can be written as
\begin{eqnarray}
\min_{|T|\leq s, T\subseteq[m]}\Big\{\min_{\bfal\in \R^{m}}   D(\alpha):\langle\bfal_T, \bfy_T \rangle=0\Big\}.
\end{eqnarray}
It follows from \eqref{hessian} that $D(\cdot)$ is strongly convex.  Therefore, the inner problem is a strongly convex program that admits a unique solution, say $\alpha_T$. In addition, the choices of $T$ such that $|T|\leq s, T\subseteq[m]$ are finitely many. To derive the global optimal solution, we pick one $T$ from the choices making $D(\alpha_T)$ the smallest.  \hfill   $ \Box $

\subsection*{Proof of  \cref{stationary-equation}}

{\bf proof}    It follows from $\bfz^*$ being an $\eta$-stationary point and \eqref{SM-SVM-h-equ-dual-sparse-tau} that $\langle\bfal^*, \bfy \rangle =0$ and 
\begin{eqnarray*} \bfal^*& \in& \P_s\left(\bfal^*-\eta g ( \bfz^*)\right)\nonumber\\
&\overset{\eqref{P-s-T}}{=}&\left\{\left[\begin{array}{c}
\bfal^*_T-\eta g_T ( \bfz^*)\\
0
\end{array}\right]: T\in \T_s(\bfal^*-\eta g ( \bfz^*))\right\},
\end{eqnarray*}
which is equivalent to stating that there exists a $T_*\in \T_s(\bfal^*-\eta g ( \bfz^*)  )$ satisfying $\bfal_{\overline T_*}^*=0$ and $0= g_{  T_*} ( \bfz^*).$
This can reach  the conclusion immediately. \hfill   $ \Box $
\subsection*{A Lemma for \cref{nec-suff-opt-con-tau} }
To prove \cref{nec-suff-opt-con-tau}, we need the following Lemma.
\begin{lemma}\label{nec-suff-opt-con-KKT} Consider a feasible point $\bfal^*$ to the problem \eqref{SM-SVM-h-equ-dual-sparse}, namely, $\|\bfal^*\|_0\leq s$ and $\langle\bfal^*, \bfy \rangle=0$.
\begin{itemize}
\item[$a)$] It is a local minimizer if and only if there is  $ b ^*\in\R$ such that \begin{eqnarray}
\label{SM-SVM-h-equ-dual-sparse-opt-kkt}
\left\{\begin{array}{rll}
 g_{S_*}(\bfz^*) &=&0, \\
 \langle\bfal^*, \bfy \rangle&=&0, \\
  \|\bfal^*\|_0&=& s,  
\end{array}\right.
~~{\rm or}~~
\left\{\begin{array}{rll}
g (\bfz^*)&=&0,  \\
 \langle\bfal^*, \bfy \rangle&=&0, \\
  \|\bfal^*\|_0&<& s. 
\end{array}\right.
\end{eqnarray}
\item[$b)$] If $\|\bfal^*\|_0<s$, then the local minimizer and the global minimizer are identical to each other and unique.
\end{itemize}
\end{lemma}
{\bf proof} We only prove that the global  minimizer $\bfal^*$ is unique if $\|\bfal^*\|_0<s$. The proofs of the remaining parts can be seen in \cite[Theorem 3.2]{PFX17}. The  condition in \eqref{hessian} indicates $D(\bfal)$ is a strongly convex function and thus enjoys the property
\begin{eqnarray}\label{strong-convex}
D(\bfal)%&=& D(\bfal') + \langle \nabla D(\bfal'), \bfal-\bfal'\rangle\nonumber\\
%& +& \langle \bfal-\bfal', H(\bfal_t) (\bfal-\bfal')\rangle/2,\nonumber\\
 &\overset{\eqref{hessian}}{\geq}&  D(\bfal') + \langle \nabla D(\bfal'), \bfal-\bfal'\rangle \nonumber\\
 &+& \langle \bfal-\bfal', P(\bfal-\bfal')\rangle/2,
\end{eqnarray}
for any $\bfal,\bfal'\in\R^m$. If there is another global  minimizer  $\bfal\neq \bfal^*$, then the strong convexity  of $ D(\cdot)$  gives  rise to
\begin{eqnarray*}
&&D(\bfal)- D(\bfal^*)\\
 &\overset{\eqref{strong-convex}}{\geq}&   \langle \bfal-\bfal^*,P(\bfal-\bfal^*)\rangle/2  + \langle \nabla D(\bfal^*), \bfal-\bfal^*\rangle  \\
%&=&  D(\bfal^*) + \langle\bfal-\bfal^*,P(\bfal-\bfal^*)\rangle/2+ \langle P   \bfal^*-{\bf1}, \bfal-\bfal^*\rangle  \\
&\overset{\eqref{beta-mu}}{=}&    \langle \bfal-\bfal^*,P(\bfal-\bfal^*)\rangle/2  +\langle  g( \bfz^*)-\bfy b ^*, \bfal-\bfal^*\rangle  \\
&\overset{\eqref{SM-SVM-h-equ-dual-sparse-opt-kkt}}{=}&    \langle \bfal-\bfal^*,P(\bfal-\bfal^*)\rangle/2  - \langle  \bfy b ^*, \bfal-\bfal^*\rangle  \\
&\overset{\eqref{SM-SVM-h-equ-dual-sparse-opt-kkt}}{=}&   \langle \bfal-\bfal^*,P(\bfal-\bfal^*)\rangle/2\\
&  \overset{\eqref{hessian}}{>}&    0,
\end{eqnarray*}
where the third equation is valid because global  minimizers $\bfal$ and $\bfal^*$  satisfy $\langle\bfal^*, \bfy \rangle=\langle\bfal, \bfy \rangle=0$. It follows from the global optimality that $D(\bfal^*) = D(\bfal)$, contradicting  the above inequality. Therefore, $\bfal^*$ is unique.   \hfill   $ \Box $

\subsection*{Proof of  \cref{nec-suff-opt-con-tau} }
{\bf proof} a) It follows from \cite[Lemma 2.2]{BE13} that an $\eta$-stationary point in \eqref{SM-SVM-h-equ-dual-sparse-tau}  can be equivalently written as
\begin{eqnarray}
\label{SM-SVM-h-equ-dual-sparse-tau-equ}
\left\{\begin{array}{rlr}
g_{S_*}( \bfz^*) &=&0,\\
 \eta \|g_{\overline S_*} ( \bfz^*) \|_{[1]}&\leq& \|\bfal^*\|_{[s]}, \\
  \|\bfal^*\|_0&\leq& s, \\
   \langle\bfal^*, \bfy \rangle&=&0.
  \end{array}\right.
\end{eqnarray}
 This  clearly indicates \eqref{SM-SVM-h-equ-dual-sparse-opt-kkt} by $\|\bfal^*\|_{[s]}=0$ if $\|\bfal^*\|_0<s$. Therefore it is a local minimizer by  \cref{nec-suff-opt-con-KKT} a).

b) If $\|\bfal^*\|_0<s$, then  \cref{nec-suff-opt-con-KKT} a) states that a local minimizer  $\bfal^*$  satisfies the second condition in \eqref{SM-SVM-h-equ-dual-sparse-opt-kkt} which is same as \eqref{SM-SVM-h-equ-dual-sparse-tau-equ}. Therefore, it is also an $\eta$-stationary point for any $\eta>0$. If $\|\bfal^*\|_0=s$, then $\|\bfal^*\|_{[s]}>0$, which implies 
$$\eta^*:=  {\|\bfal^*\|_{[s]}}/({ 2\|H( \bfal^*) \bfal^*-\bfe\|_{[1]} })>0.$$
 A local minimizer satisfies the first condition in \eqref{SM-SVM-h-equ-dual-sparse-opt-kkt},  
\begin{eqnarray*}g_{S_*}( \bfz^*)\overset{\eqref{beta-mu} }{=}(H( \bfal^*) \bfal^*)_{S_*}-\bfe+\bfy_{S_*} b ^*=0,\end{eqnarray*}
which gives rise to
$$| b ^*|=\|\bfy_{  S_*}\|_{[1]}| b ^*|=\| (H( \bfal^*) \bfal^*)_{S_*}-\bfe \|_{[1]}$$
because of $|\bfy|={\bf1}$. In addition, $0<\eta<\eta^*$ gives rise to 
\allowdisplaybreaks\begin{eqnarray*}
&&\|g_{\overline S_*}( \bfz^*)\|_{[1]}\\
 &\overset{\eqref{beta-mu} }{=}&
  \| (H( \bfal^*) \bfal^*)_{\overline S_*}-\bfe  +\bfy_{\overline S_*} b ^*\|_{[1]}\\
 &\leq&  \|  (H( \bfal^*) \bfal^*)_{\overline S_*}-\bfe  \|_{[1]}+| b ^*|\|\bfy_{\overline S_*}\|_{[1]}  \\
  &=& \|(H( \bfal^*)\bfal^*)_{\overline S_*}-\bfe  \|_{[1]}+| b ^* |  \\
    &=&  \| (H( \bfal^*)\bfal^*)_{\overline S_*}-\bfe  \|_{[1]}+\|  
    (H( \bfal^*)\bfal^*)_{  S_*}-\bfe  \|_{[1]}  \\
 &\leq& \|  H( \bfal^*)\bfal^*-\bfe  \|_{[1]} 2\\
 & =& \|\bfal^*\|_{[s]}/\eta^* \\
  &\leq&  \|\bfal^*\|_{[s]}/\eta.
\end{eqnarray*}
This verifies the second inequality in \eqref{SM-SVM-h-equ-dual-sparse-tau-equ}, together with \eqref{SM-SVM-h-equ-dual-sparse-opt-kkt}  claiming the conclusion. 

 c) An $\eta$-stationary point $\bfal^*$ with $\|\bfal^*\|_0<s$ satisfies the condition \eqref{SM-SVM-h-equ-dual-sparse-tau-equ}, which is same as the second case $\|\bfal^*\|_0<s$ in \eqref{SM-SVM-h-equ-dual-sparse-opt-kkt}. Namely, $\bfal^*$ is also a local minimizer, which makes the conclusion immediately from   \cref{nec-suff-opt-con-KKT} b). 
 
 d) An $\eta$-stationary point $\bfal^*$ with $\|\bfal^*\|_0=s$ satisfies 
 \begin{eqnarray*} 
 ~~~~\bfal^*\overset{\eqref{SM-SVM-h-equ-dual-sparse-tau}}{\in} \P_s(\bfal^*-\eta g ( \bfz^*)) \overset{\eqref{HTP-operator}}{=} \underset{ \|\bfal\|_0\leq s}{\rm argmin}\|\bfal- (\bfal^*-\eta g ( \bfz^*) )\| , \end{eqnarray*}
which means for any feasible point, namely, $\|\bfal\|_0\leq s$ and $\langle\bfal, \bfy \rangle=0$, we have
 \begin{eqnarray*}
\left\|\bfal^*- (\bfal^*-\eta g ( \bfz^*) )\right\|^2 \leq \left\|\bfal - (\bfal^*-\eta g ( \bfz^*) )\right\|^2.
\end{eqnarray*}
This suffices to 
 \begin{eqnarray}
\label{fact1}
-\|\bfal^*-\bfal\|^2 \leq  2\eta \langle \bfal- \bfal^*,   g ( \bfz^*)\rangle.
\end{eqnarray} 
The strong  and quadratic convexity of $ D(\cdot)$   gives  rise to
\begin{eqnarray*}
&&2D(\bfal)-2D(\bfal^*) \\
&\overset{\eqref{strong-convex}}{\geq}&  \langle\bfal-\bfal^*,P(\bfal-\bfal^*)\rangle + 2\langle \nabla D(\bfal^*), \bfal-\bfal^*\rangle  \\
&\overset{\eqref{beta-mu}}{=}&  \langle\bfal-\bfal^*,P(\bfal-\bfal^*)\rangle + 2\langle  g ( \bfz^*)-\bfy b ^*, \bfal-\bfal^*\rangle  \\
&\overset{\eqref{SM-SVM-h-equ-dual-sparse-tau-equ}}{=}& \langle\bfal-\bfal^*,P(\bfal-\bfal^*)\rangle + 2\langle  g ( \bfz^*), \bfal-\bfal^*\rangle  \\
&\overset{\eqref{fact1}}{\geq}& \langle\bfal-\bfal^*,P(\bfal-\bfal^*)\rangle -\|\bfal^*-\bfal\|^2/ \eta   \\
&=& \langle\bfal-\bfal^*,(P-I/\eta)(\bfal-\bfal^*)\rangle\\
& \geq&  0,
\end{eqnarray*}
where the last inequity follows from 
 $$P-{I}/{\eta} =\left[{1}/{ C}-{1}/{\eta}\right]I+Q^\top Q\succeq 0$$
by $\eta\geq C$. Therefore, $\bfal^*$ is a global minimizer. If  there is another global minimizer $\hat\bfal\neq \bfal^*$, then the strictness $\succ$ in the above condition by $\eta> C$ leads to a contradiction,
 \begin{eqnarray*} 0 &=& D(\hat\bfal)-D(\bfal^*)\\
 & \geq&  \langle\hat\bfal-\bfal^*,(P-I/\eta)(\hat\bfal-\bfal^*)\rangle\\
 &>&0.\end{eqnarray*}
Hence, $\bfal^*$ is unique. The  proof is completed. \hfill   $ \Box $

\subsection*{A Lemma for \cref{the:i-iterate-convergence}}
Before proving \cref{the:i-iterate-convergence}, 
%we define some constants:
%\allowdisplaybreaks\begin{eqnarray}
% \gamma~&:=&2\max\Big\{ 1+\eta/c+\eta\|Q \|^2,~ \eta\sqrt{m}\Big\},\nonumber\\
%\label{tau*}\eta^*&:=&\begin{cases}
%  \|\bfal^*\|_{[s]} \|g(\bfz^*)\|_{[1]}^{-1} , &\| \bfal^*\|_0=s,\\
%+\infty, & \| \bfal^*\|_0<s,
%\end{cases} \\
% \delta^*&:=&\begin{cases}
%\gamma^{-1} ( \|\bfal^*\|_{[s]}- \eta\|g(\bfz^*)\|_{[1]} ), & \| \bfal^*\|_0=s,\\
% \gamma^{-1}  \min_{i\in S_*}|\alpha_i^*|, & 0<\| \bfal^*\|_0<s,\\
%+\infty, & \| \bfal^*\|_0 =0.
%\end{cases} \nonumber
%\end{eqnarray}
%Based on which, 
we first present some properties regarding an $\eta$-stationary point of \eqref{SM-SVM-h-equ-dual-sparse}.
\begin{lemma}\label{lemma-31} Let $\bfz^*$ be an $\eta$-stationary point of \eqref{SM-SVM-h-equ-dual-sparse} with $0<\eta<\eta^*$, where $\eta^*$   is given by \eqref{tau*}.  Then there  always exists a $\delta^*>0$ such that for any $\bfz\in N(\bfz^*,\delta^*)$ with $\|\bfal\|_0\leq s$,  the following results hold.
\begin{itemize}
\item[$a)$] If $\| \bfal^*\|_0=s$, then 
\begin{eqnarray}
\label{relation-S-T}
 \T_s(\bfal-\eta g(\bfz) = \T_s(\bfal^*-\eta g(\bfz^*) )=\{ S_*  \}.
\end{eqnarray}
 If $\| \bfal^*\|_0<s$, then for any $T\in \T_s(\bfal-\eta g(\bfz) )$ and any ${T_*}\in \T_s(\bfal^*-\eta g(\bfz^*) )$, it holds
\begin{eqnarray}
\label{relation-S-T-1}
S_* \subseteq (T_*\cap T).
\end{eqnarray}
\item[$b)$] For any $T\in \T_s(\bfal-\eta g(\bfz) )$, it holds
\begin{eqnarray}\label{stationary-equation-T-1}
F(\bfz^*;T)=0.
\end{eqnarray}
\end{itemize}
\end{lemma}
{\bf proof}  a)   It follows from  \cref{stationary-equation} and $\bfz^*$ being an $\eta$-stationary point of \eqref{SM-SVM-h-equ-dual-sparse} that  $F(\bfz^*;T_*)=0$  
%\begin{eqnarray}
%\label{stationary-equation-T-*}
%F(\bfz^*;T_*)=\left[
%\begin{array}{c}
%  g_{T_*}(\bfz^*) \\ 
%  \bfal_{\overline T_*}^*  \\
%   \langle\bfal_{T_*}^*, \bfy_{T_*} \rangle
%\end{array}\right]=0.
%\end{eqnarray}
for any ${T_*}\in \T_s(\bfal^*-\eta g(\bfz^*) )$. We first derive  
\begin{eqnarray}\label{EE0}
E(\bfal^*)\bfal^*=E(\bfal)\bfal^*.
\end{eqnarray}
If $\bfal^*=0$, it is true clearly. If $\bfal^*\neq 0$, we have 
 \begin{eqnarray}\label{EE00}
\alpha_i^*>0 \Longrightarrow \alpha_i>0,~~~~
\alpha_i^*<0 \Longrightarrow \alpha_i<0.
\end{eqnarray}
If \eqref{EE00} is not true, then there is a  $j$ that violates one of the above relations, namely $\alpha_j^*$ and $\alpha_j$ have different signs. As a consequence,  
$$ |\alpha_j^*|<|\alpha_j^* -\alpha_j|\leq \|\bfal^*-\bfal\|<\delta^*.$$
This is a contradiction for a sufficiently small positive $\delta^*$.
Therefore, we must have \eqref{EE00}. Recall that $E(\bfal)$ is a diagonal matrix with diagonal elements given by \eqref{diag-hessian}. It follows
\begin{eqnarray*} 
&& [(E(\bfal^*)-E(\bfal))\bfal^* ]_i\\ 
& =&  [E_{ii}(\bfal^*)-E_{ii}(\bfal) ]\alpha_i^*\\
&\overset{\eqref{EE00},\eqref{diag-hessian}}{=}&\left\{
\begin{array}{rr}
(1/C-1/C)\alpha_i^*,&\alpha_i^*>0 \\
(1/c-1/c)\alpha_i^*,&\alpha_i^*<0 \\
(E_{ii}(\bfal^*)-E_{ii}(\bfal))0,&\alpha_i^*=0 \\
\end{array}
\right.\\ 
& =& 0.
\end{eqnarray*}
Therefore, \eqref{EE0} is true, allowing us to derive that
  \begin{eqnarray} 
&& \|H(\bfal^*)\bfal^*-H(\bfal)\bfal\|\nonumber\\
&\overset{\eqref{gradient-hessian}}{=}&  \|(E(\bfal^*)+Q^\top Q )\bfal^*-(E(\bfal)+Q^\top Q )\bfal \|\nonumber\\
&\leq& \|E(\bfal^*)\bfal^*-E(\bfal)\bfal\|+ \|Q \|^2 \|\bfal^*-\bfal\|\nonumber\\
&\overset{\eqref{EE0}}{=}& \|E(\bfal)(\bfal^*-\bfal)\|+ \|Q \|^2 \|\bfal^*-\bfal\|\nonumber\\ 
&\overset{\eqref{diag-hessian}}{\leq}&  (1/c+\|Q \|^2)\|\bfal^*-\bfal\|\nonumber\\ 
\label{contrad-s-00}& \leq&  (1/c+\|Q \|^2) \delta^*, 
\end{eqnarray} 
where $\|Q \|$ is the spectral norm of $Q$. This suffices to
{\allowdisplaybreaks \begin{eqnarray} 
&&\|g(\bfz^*)-g(\bfz)\|\nonumber\\
&\overset{\eqref{beta-mu}}{=}& \|H(\bfal^*)\bfal^*-H(\bfal)\bfal+( b ^*- b )\bfy \|\nonumber\\ 
&\leq& \|H(\bfal^*)\bfal^*-H(\bfal)\bfal\|+| b ^*- b |\|\bfy \|\nonumber\\ 
\label{contrad-s-0}&\overset{\eqref{contrad-s-00}}{\leq}& (1/c+\|Q \|^2+\sqrt{m}) \delta^*. 
\end{eqnarray}}
{\bf Case i)} $\|\bfal^*\|_0=s$. 
Since $|T_*|=s$ and $\bfal_{\overline T_*}^*=0$ from \eqref{stationary-equation-T}, it holds $T_*=\supp(\bfal^*)=S_*$.
If $\|  g_{\overline S_*}(\bfz^*)\|_{[1]}=0$, then 
\begin{eqnarray}
\label{T*-S*-0}
 \|\bfal^*_{\overline S_*}-\eta g_{\overline S_*}(\bfz^*)\|_{[1]}&\overset{\eqref{stationary-equation-T}}{=}&0<
\|\bfal^*_{S_*}\|_{[s]}\nonumber\\
&\overset{\eqref{stationary-equation-T}}{=}&\|\bfal^*_{S_*}-\eta g_{S_*}(\bfz^*)\|_{[s]}
\end{eqnarray}
If $\|  g_{\overline S_*}(\bfz^*)\|_{[1]}\neq0$, then
\begin{eqnarray*}
 \|\bfal^*_{\overline S_*}-\eta g_{\overline S_*}(\bfz^*)\|_{[1]}&\overset{\eqref{stationary-equation-T}}{=}&\eta\|  g_{\overline S_*}(\bfz^*)\|_{[1]} \\
 & <&  \eta^*\|  g_{\overline S_*}(\bfz^*)\|_{[1]} \\
 &\overset{\eqref{tau*}}{=}&\|\bfal^*_{S_*}\|_{[s]}\\
   &\overset{\eqref{stationary-equation-T}}{=}& 
\|\bfal^*_{S_*}-\eta g_{S_*}(\bfz^*)\|_{[s]}
\end{eqnarray*}
Therefore, both scenarios derive \eqref{T*-S*-0}. This means $S_*$ contains the $s$ largest elements of $|\bfal^*_{S_*}-\eta g_{S_*}(\bfz^*)|$, 
together with the definition of  
$\T_s$ in \eqref{T-beta} indicating $$\T_s(\bfal^*-\eta g(\bfz^*) )= \{ S_*  \}.$$
%Next, we check $S_*\subseteq \supp(\bfal)$. In fact, if it is not true, then there is an $i\in S_*$ but $i\notin \supp(\bfal)$, which  incurs  
%\begin{eqnarray}
%\label{ S*-supp-alpha} \|\bfal^*\|_{[s]}  \leq  |\alpha_i^*|=|\alpha_i^*-\alpha_i| \leq \|\bfal^*-\bfal\| <\delta^*.\end{eqnarray}
%This is a contradiction for a sufficiently small positive $\delta^*$. So  $S_*\subseteq \supp(\bfal)$, together with  $|S_*|=s \geq \|\bfal\|_0=| \supp(\bfal)|$ yielding  $S_*= \supp(\bfal)$.  
Next, we show $\T_s(\bfal-\eta g(\bfz) )=  \{ S_*  \}$ , i.e., to show \begin{eqnarray}
\label{T*-S*-00}
 \|\bfal_{\overline S_*}-\eta g_{\overline S_*}(\bfz)\|_{[1]}<
\|\bfal_{S_*}-\eta g_{S_*}(\bfz)\|_{[s]}
\end{eqnarray}
 For sufficiently small $\delta^*$, $ \bfal$ can be close enough to $ \bfal^* $ and $g(\bfz)$ can be close enough to $g(\bfz^*)$ from \eqref{contrad-s-0}. These and \eqref{T*-S*-00} guarantee \eqref{T*-S*-00} immediately.

{\bf Case ii)} $\|\bfal^*\|_0<s$. The fact $\bfal_{\overline T_*}^*=0$ from \eqref{stationary-equation-T}  indicates $S_*\subseteq T_*$. We next show $S_*\subseteq T$. If $\bfal^*=0$, then $S_*=\emptyset\subseteq T$ clearly. Therefore, we focus on $\bfal^*\neq0$.  
 Since $\|\bfal^*\|_0<s$, $\|\bfal^*\|_{[s]}=0$, which together with \eqref{SM-SVM-h-equ-dual-sparse-tau-equ} derives 
\begin{eqnarray}\label{beta-0}
 g (\bfz^*)=0,
\end{eqnarray}
implying that the indices of $s$ largest elements of $|\bfal^*-\eta g(\bfz^*)|=|\bfal^*|$ contain $S_*$. Note that
$$\bfal-\eta g(\bfz)=(\bfal_{S_*}-\eta g_{S_*}(\bfz);\bfal_{\overline S_*}-\eta g_{\overline S_*}(\bfz) )$$
Again, for sufficiently small $\delta^*$, $ \bfal$ and $g(\bfz)$ can be close enough to $ \bfal^* $  and $g(\bfz^*)$ from \eqref{contrad-s-0}. Therefore, $\bfal_{S_*}-\eta g_{S_*}(\bfz)$ and $\bfal_{\overline S_*}-\eta g_{\overline S_*}(\bfz)$ are close to $\bfal^*_{  S_*}-\eta g_{ S_*}(\bfz^*) =\bfal^*_{  S_*}$ and $\bfal^*_{\overline S_*}-\eta g_{\overline S_*}(\bfz^*) =0$, respectively. These imply $S_* \subseteq T$.

%If $S_* \nsubseteq T$, then the fact $|S_*|< s=|T|$ also indicates that there is an $i\in S_*, i\notin T$ and a $j\notin S_*, j\in T$. This together with the definition \eqref{T-beta} of $\T_s$ results in
%\begin{eqnarray}\label{alpha-beta-0}
%|\alpha_j-\eta  g_j(\bfz )| &\geq& |\alpha_i-\eta  g_i(\bfz )|,\nonumber \\
%|\alpha_i^*|\overset{\eqref{beta-0}}{=}|\alpha_i^*-\eta  g_i(\bfz^*)| &>&0 \overset{\eqref{beta-0}}{=} |\alpha_j^*-\eta  g_j(\bfz^*)|,
%\end{eqnarray}
%which leads to the following contradiction
%\begin{eqnarray*} 
%0< |\alpha_i^*| 
%&\overset{\eqref{alpha-beta-0}}{=}&|\alpha_i^*-\eta  g_i(\bfz^*)| - |\alpha_j^*-\eta  g_j(\bfz^*)|\\
%&\overset{\eqref{alpha-beta-0}}{\leq}&|\alpha_i^*-\eta g_i(\bfz^*)| - |\alpha_i-\eta g_i(\bfz)|\nonumber\\
%&+& |\alpha_j-\eta g_j(\bfz)|- |\alpha_j^*-\eta g_j(\bfz^*)|  \nonumber\\
%&\leq& |\alpha_i^*-\eta g_i(\bfz^*) -(\alpha_i-\eta g_i(\bfz)) |\nonumber\\ &+&  |\alpha_j^*-\eta g_j(\bfz^*) -(\alpha_j-\eta g_j(\bfz) ) |\nonumber\\
%&\leq&|\alpha_i^*-\alpha_i|+\eta|g_i(\bfz^*)-g_i(\bfz)|\\
%&+&|\alpha_j^*-\alpha_j|+\eta|g_j(\bfz^*)-g_j(\bfz)|.
%\end{eqnarray*}
%This is a contradiction since for sufficiently small $\delta^*$, $ \bfal$ and $g(\bfz)$  can be close enough to $ \bfal^* $  and $g(\bfz^*)$ from \eqref{contrad-s-0}.

b) To prove $F(\bfz^*;T)=0$, we need to show
\begin{eqnarray}
\label{stationary-equation-*-T}
F(\bfz^*;T)=\left[
\begin{array}{c}
   g _{T}(\bfz^*)  \\ 
  \bfal_{\overline T}^*  \\
   \langle\bfal_{T}^*, \bfy_{T} \rangle
\end{array}\right]=0.
\end{eqnarray}
If $\|\bfal^*\|_0=s$, then  $T=S_*=T_*$ by a), deriving the desired result by \eqref{stationary-equation-T} immediately.  If $\|\bfal^*\|_0<s$, then $ g _{T}(\bfz^*)=0$ by \eqref{beta-0}. Again from b), $S_*\subseteq (T\cap T_*)$ means $\overline T \subseteq  \overline S_*$,   indicating that $\bfal_{\overline T}^*=0$ due to  $\bfal_{\overline S_*}^*=0$. Finally, 
\begin{eqnarray*}\langle\bfal_{T}^*, \bfy_{T} \rangle&=&\langle\bfal^*_{S_*}, \bfy_{S_*} \rangle+\langle\bfal^*_{T\setminus S_*}, \bfy_{T\setminus S_*}\rangle=\langle\bfal^*_{S_*}, \bfy_{S_*} \rangle\\
&=&  \langle\bfal^*_{T_*}, \bfy_{T_*}\rangle-\langle\bfal^*_{T_*\setminus S_*}, \bfy_{T_*\setminus S_*}\rangle \overset{\eqref{stationary-equation-T}}{=} 0.\end{eqnarray*}
The proof is finished.  
\hfill   $ \Box $

\subsection*{Proof of  \cref{the:i-iterate-convergence}}
{\bf proof}  
Consider a point $\bfz_t^k:=\bfz^*+t(\bfz^k-\bfz^*)$ with $t\in[0,1]$. Since  $\bfz^{k}\in N(\bfz^*,\delta^*)$, it also holds $\bfz_t^{k}\in N(\bfz^*,\delta^*)$ due to
$$\|\bfz_t^k-\bfz^*\|=t\|\bfz^k-\bfz^*\|\leq \|\bfz^k-\bfz^*\|<\delta^*.$$
We first prove that
\begin{eqnarray}
\label{Ea-Ea}
E(\bfal^k)= E(\bfal_t^k).
\end{eqnarray}
In fact, if  $\bfal^*=0$, then $\bfal_t^k=t\bfal^k$, implying that $\bfal_t^k$ and $\bfal^k$ have the same signs. This together with the definition \eqref{diag-hessian} of $E(\cdot)$ shows \eqref{Ea-Ea} immediately. If $\bfal^*\neq0$, then the same reasoning proving \eqref{EE00} also derives that
 \begin{eqnarray*} 
\alpha_i^*>0 &\Longrightarrow& \alpha_i^k>0,~~  (\alpha_t^k)_i=(1-t)\alpha_i^*+t\alpha_i^k>0,\\
\alpha_i^*<0 &\Longrightarrow& \alpha_i^k<0,~~  (\alpha_t^k)_i=(1-t)\alpha_i^*+t\alpha_i^k<0,\\
\alpha_i^*=0 &\Longrightarrow&  ~~~~~~~~~~~~~~~(\alpha_t^k)_i= t\alpha_i^k. 
\end{eqnarray*}
These also mean  $\bfal_t^k$ and $\bfal^k$ have the same signs. Namely, \eqref{Ea-Ea} is true and results in
\begin{eqnarray*}
~~H(\bfal^k) &\overset{\eqref{gradient-hessian}}{=}& E(\bfal^k)+Q^\top Q\\
& \overset{\eqref{Ea-Ea}}{=}&  E(\bfal_t^k)+Q^\top Q=H(\bfal_t^k),
\end{eqnarray*}
contributing to $
\nabla  F (\bfz^{k}_t;T_k)=\nabla  F (\bfz^{k};T_k)$ due to 
\begin{eqnarray*}
\left[
\begin{array}{ccc}
 H_{ T_k}(\bfal_t^k) & 0 & \bfy_{T_k}\\ 
0&I&0\\
  \bfy_{T_k}^\top&0&0  
\end{array}\right] 
=\left[
\begin{array}{ccc}
 H_{ T_k}(\bfal^k) & 0 & \bfy_{T_k}\\ 
0&I&0\\
  \bfy_{T_k}^\top&0&0  
\end{array}\right].
\end{eqnarray*}
for any $T_k\in \T_s(\bfal^k-\eta g (\bfz^k) )$.  
It follows from the mean value theorem that there exists a $\bfz^{k}_t$ satisfying
 \begin{eqnarray*} F (\bfz^{k};T_k)&\overset{\eqref{stationary-equation-T-1}}{=}&F (\bfz^{k};T_k) -F (\bfz^{*};T_k)\\
 &=& \nabla  F (\bfz^{k}_t;T_k) (\bfz^{k}-\bfz)\\
 &=& \nabla  F (\bfz^{k};T_k) (\bfz^{k}-\bfz^{*}). 
\end{eqnarray*}
This together with $\nabla F (\bfz^{k};T_k)$ being always nonsingular due to \eqref{Jacobian-F} enables to show that
\begin{eqnarray*}
\bfz^* &=& \bfz^k-(\nabla F (\bfz^{k};T_k))^{-1}  F (\bfz^{k};T_k)\\ &\overset{\eqref{Newton-Method}}{=}&  \bfz^k +\bfd^k \overset{\eqref{form4}}{=}  \bfz^{k+1}.
\end{eqnarray*}
Finally, for any $T_{k+1}\in \T_s(\bfal^{k+1}-\eta g  (\bfz^{k+1}) )= \T_s(\bfal^{*}-\eta g   (\bfz^*) )$, it follows from  $\bfz^*$ being an $\eta$-stationary point that
$$\|F(\bfz^{k+1},T_{k+1})\|=\|F(\bfz^*,T_{k+1})\|\overset{\eqref{stationary-equation-T}}{=}0.$$
The entire proof is completed.  
\hfill   $ \Box $

\end{document}